\documentstyle[11pt, bezier]{article}
\newcounter{ppp}
\newcounter{pdfour}
\newcounter{pdfive}

\newcounter{pdsix}
\newcounter{pdseven}
\newcounter{pdeight}
\newcounter{pdnine}
\newcounter{pdten}
\newcounter{pdeleven}
\newcounter{pdtwelve}
\newcounter{pdthirteen}
\newcounter{pdfourteen}

\newcounter{pdfifteen}
\newcounter{pdsixteen}
\newcounter{pdseventeen}

\newcommand{\Lab}{{\rm Lab}}

\newcommand{\ttt}{{\cal T}}

\newcommand{\aaa}{{\cal A}}

\newcommand{\bott}{{\bf bot}}
\newcommand{\vk}{van Kampen }

\newcommand{\iv}{^{-1}}

\newcommand{\rr}{{\cal R} }

\newcommand{\qq}{{\cal Q} }
\newcommand{\sss}{{\cal S} }

\begin{document}

\newtheorem{lemma}{Lemma}
\newtheorem{theorem}{Theorem}
\newtheorem{cor}{Corollary}
\newtheorem{definition}{Definition}

\title{Length and Area Functions on Groups and Quasi-Isometric Higman
Embeddings}
\author{A.Yu. Olshanskii and M.V. Sapir\thanks{The research of
the first author was supported in part by the Russian fund for fundamental
research 96-01-420. The research of the second author was supported in part
by the NSF grant DMS 9623284}}

\date{}
\maketitle

\begin{abstract}
We survey recent results about asymptotic functions of groups, obtained by
the authors in collaboration with J.-C.Birget, V. Guba and E. Rips. We also
discuss methods used in the proofs of these results.
\end{abstract}
\tableofcontents
\section{Results}

\subsection{Definitions}

Recall that isoperimetric functions of a finitely presented group $G=\langle
X\ | \ R\rangle$ measure areas of van Kampen diagrams over the presentation of
this group.  \setcounter{pdfour}{\value{ppp}} Figure \thepdfour\ shows what
a van Kampen diagram may look like.

\bigskip

\unitlength=1.00mm
\linethickness{0.4pt}
\begin{picture}(118.47,59.67)
\put(80.17,29.83){\oval(75.67,51.00)[]}
\put(42.33,17.33){\line(1,0){30.00}}
\put(72.33,17.33){\line(0,-1){13.00}}
\put(72.33,4.33){\line(0,1){6.00}}
\put(72.33,10.33){\line(1,0){22.33}}
\put(94.67,10.33){\line(0,-1){6.00}}
\put(94.67,4.33){\line(0,1){19.00}}
\put(94.67,23.33){\line(1,0){23.33}}
\put(118.00,23.33){\line(-1,0){45.67}}
\put(72.33,23.33){\line(0,-1){6.00}}
\put(72.33,17.33){\line(0,1){18.00}}
\put(72.33,35.33){\line(-1,0){30.00}}
\put(42.33,35.33){\line(1,0){16.00}}
\put(58.33,35.33){\line(0,1){20.00}}
\put(58.33,55.33){\line(0,-1){10.00}}
\put(58.33,45.33){\line(1,0){26.33}}
\put(84.67,45.33){\line(0,1){10.00}}
\put(84.67,55.33){\line(0,-1){20.00}}
\put(84.67,35.33){\line(-1,0){12.33}}
\put(72.33,35.33){\line(1,0){36.33}}
\put(108.67,35.33){\line(0,-1){18.33}}
\put(108.67,17.00){\line(-1,0){14.00}}
\put(94.67,17.00){\line(0,1){18.33}}
\put(100.33,35.33){\line(0,1){20.00}}
\put(73.67,59.67){\makebox(0,0)[cc]{$w$}}
\put(91.67,44.67){\makebox(0,0)[cc]{$r_i$}}
\put(51.33,35.33){\circle*{0.67}}
\put(58.33,41.00){\circle*{0.67}}
\put(66.00,45.33){\circle*{0.67}}
\put(75.33,45.33){\circle*{0.67}}
\put(72.33,29.00){\circle*{0.67}}
\put(84.00,23.33){\circle*{0.67}}
\put(116.33,53.33){\circle*{0.94}}
\put(43.67,52.67){\circle*{0.67}}
\put(56.33,17.33){\circle*{0.67}}
\put(101.67,23.33){\circle*{0.67}}
\put(108.67,23.33){\circle*{0.67}}
\put(116.00,6.33){\circle*{0.94}}
\put(44.00,6.33){\circle*{0.67}}
\put(42.33,26.33){\circle*{0.67}}
\put(118.00,36.33){\circle*{0.94}}
\put(45.67,33.67){\makebox(0,0)[cc]{$a$}}
\put(55.00,33.67){\makebox(0,0)[cc]{$b$}}
\put(43.67,42.33){\makebox(0,0)[cc]{$b$}}
\put(51.33,53.67){\makebox(0,0)[cc]{$a$}}
\put(57.00,38.33){\makebox(0,0)[cc]{$c$}}
\put(57.33,49.67){\makebox(0,0)[cc]{$d$}}
\put(57.00,43.33){\makebox(0,0)[cc]{$d$}}
\put(42.33,41.00){\vector(0,1){1.33}}
\put(42.33,42.33){\vector(0,0){0.00}}
\put(45.00,35.33){\line(1,0){3.67}}
\put(51.67,55.33){\line(1,0){4.67}}
\put(58.33,48.33){\line(0,1){3.33}}
\put(51.67,55.33){\vector(1,0){3.67}}
\put(58.33,48.67){\vector(0,1){4.67}}
\put(58.33,44.67){\vector(0,-1){2.33}}
\put(58.33,39.00){\vector(0,-1){2.67}}
\put(56.33,35.33){\vector(-1,0){3.33}}
\put(45.00,35.33){\vector(1,0){4.67}}
\end{picture}

\begin{center}
\nopagebreak[4]
Fig. \theppp.

\end{center}
\addtocounter{ppp}{1}

\bigskip
It is a directed planar labeled graph where every edge is labeled by a
generator from $X$, and the contour of every 2-cell (face) is labeled by a
relator from $R$.  By \vk lemma \cite{LS} a word $w\in X\cup X\iv$ is equal
to 1 in the group $G$ if and only if there exists a \vk diagram $\Delta$ over the
presentation of $G$ with boundary label $w$. The number of cells in $\Delta$
is equal to the number of factors in a representation of $w$ as a
product of conjugates of relators from $R$:
\begin{equation}w=\prod r_i^{u_i}.\label{eqin}\end{equation}
We are going to establish a more precise relation between equation (\ref{eqin})
and the \vk diagram $\Delta$ later (see Lemma \ref{lmm1} below).

If a van Kampen diagram has minimal number of cells, $m$, among all diagrams
with the same boundary label $w$ then we say that $w$ has {\em area} $m$. A
function $f(m)$ is called an {\em isoperimetric function} of
the presentation $\langle X \ | \ R\rangle$ of
the group $G$ if every word of length at most $m$ which is equal to 1 in the
group has area at most $f(m)$. On the set of functions
${\bf N}\to {\bf N}$, one can define a quasi-order saying that $f\prec g$
if $$f(m)\le Cg(Cm)+Cm$$ for all $m$ and some constant $C$.
Any minimal (with respect to $\prec$) isoperimetric function
is called the {\em Dehn function} of the presentation
$\langle X \ |\ R\rangle$. The article ``the" is appropriate here because
it is well known \cite{Alo},  \cite{Gersten}, \cite{gromov1} that Dehn
functions of different presentations of the same group are {\em equivalent}
that is they satisfy inequalities $$f_1\prec f_2, \hbox{\quad} f_2\prec
f_1$$.

We can also define {\em isodiametric functions} introduced by Gersten
\cite{Gersten93}. These functions measure the diameter of a van Kampen
diagram with given perimeter\footnote{Recall that the {\em diameter} of a
graph is the maximal distance between two vertices of the graph.}.
More precisely, with every word $w$ which is equal to 1 in $G$, we associate
its {\em diameter}, that is the smallest diameter of a \vk diagram with
boundary label $w$. Then if $d(n)$ is an isodiametric function of $G=\langle
X \ |\ R\rangle$, $d(n)$ must exceed the diameter every word $w= 1$ (mod
 $G$) of length $\le n$.  The equivalence of isodiametric functions is
defined as before. Isodiametric functions of different presentations of the
same group are always equivalent \cite{gromov1}.

Both isoperimetric and isodiametric functions reflect the decidability of
the word problem in the group. In particular, \cite{Gersten93}, the word
problem is decidable if and only if the Dehn function (the smallest
isodiametric function) is recursive. Nevertheless, the word problem in a
group with huge Dehn function may be easy.  For example the word problem in
the Baumslag-Solitar group $\langle a,b \ |\ a^b=a^2\rangle$ can be solved
in quadratic time (since this group is representable by $2\times 2$ integer
matrices) while the Dehn function is exponential \cite{WPG}.  One of our
main goals is to show that still there exists a very close connection
between the Dehn functions and the computational complexity of the word
problem.

If $G=\langle X\rangle$ and $H=\langle X\cup Y\rangle$ are finitely
presented, one can also define the {\em area} function  of $G$. This
function is defined on the set $W$ of all words in the alphabet $X$ which
are equal to 1 in $G$.  It takes every word from $W$ to the area of this
word in $H$.

Other  important concepts are the one of a {\em distortion function} and the
one of a {\em length function}. Let $G=\langle X \rangle$ be a finitely
generated subgroup of a finitely generated group $H=\langle Y\rangle$. Then
the distortion function $d_{G,H}(n): {\bf N}\to {\bf N}$ takes every natural
number $n$ to $\max\{|u|_{X}\ |\ u\in G, |u|_Y\le n\}$. In other words, in
order to compute $d_{G,H}(n)$ we consider all (finitely many)  elements of
$G$ whose lengths in $H$ are at most $n$, for each of these elements we
compute its length in the alphabet $X$, and then take the maximum of these
lengths.

The corresponding {\em length function} of $G$ inside $H$ is the function
$\ell: G\to {\bf N}$ which takes every element $g$ of $G$ to $|g|_Y$, the
length of $g$ in $H$.

Two functions $f_1, f_2 :  G\to {\bf N}$ are called {\em $O$-equivalent} if
$f_1(g)\le c f_2(g), f_2(g)\le cf_1(g)$ for some constant $c$ and every
$g\in G$. Different choices of generating sets in $G$ and $H$ lead to
$O$-equivalent length functions $\ell_H:G\to {\bf N}$ and equivalent
distortion functions associated with this embedding.

If $|g|_X= O(|g|_Y)$ for every $g\in G$, or, equivalently, if the distortion
function is at most linear we say that $G$ is {\em quasi-isometrically
embedded} into $H$ or that $G$ has {\em bounded distortion} in $H$.
Otherwise we say that $G$ is has unbounded distortion.

For example, every subgroup of the free group has (obviously)
bounded distortion, but the (cyclic) center $C=\langle c\rangle$ of the
3-dimensional Heisenberg group $H^3=\langle a,b,c\| \ [a,b]=c, ca=ac,
cb=bc\rangle$ has quadratic distortion. Indeed $c^{n^2}=[a^n,b^n]$ for
every $n$, the length of $c^{n^2}$ in $C$ is $n^2$ and the length of this
element in $H^3$ is $\le 4n$.

Just as Dehn functions and isodiametric functions reflect the decidability
of the word problem, the distortion function reflects
the decidability of the membership problem for subgroups: if
$G$ is a finitely generated subgroup of a finitely generated group $H$ which
has solvable word problem, then the membership in $G$ for elements of $H$ is
decidable if and only if the distortion function of $G$ in $H$ is recursive
\cite{Farb}.

In this paper, we survey recent results about isoperimetric, isodiametric,
length and area functions of groups obtained by the authors in collaboration
with J.-C. Birget, V. Guba and E. Rips.

There are several important connections between Dehn functions and length
functions. We present two easy statements without proofs here.

\begin{theorem} (Bridson, \cite{Bridson})
Let $G$ be a finitely presented group and $H\le G$ be a finitely generated
subgroup of $G$ with distortion function $d(n)$. Then the Dehn
function of the HNN extension $H = \langle G, t\ | h^t=h, h\in H\rangle$
is at least $d(n)$.
\label{disss}
\end{theorem}

\begin{theorem} (Olshanskii, Sapir 1998) The set of distortion functions of
finitely generated subgroups of the direct product of two free groups $F_2\times F_2$ coincides
(up to equivalence) with the set of all Dehn function of finitely presented
groups.
\label{mikh}
\end{theorem}

Theorem \ref{mikh} is new although a remark in \cite{gromov1} hints to a
possibility of some connection between distortion of subgroups in $F_2\times
F_2$ and Dehn functions. Here is a {\bf proof} of this theorem. It uses the well known Mikhailova's trick (see \cite{LS}) and a result from Baumslag and Roseblade \cite{BR}.

It is proved in \cite{BR} that every finitely generated subgroup $E$ of $F_2 \times F_2$
is the equalizer (in \cite{BR} it is called the free corner pullback) of two homomorphisms
$\phi: E'\to G$ and $\psi: E''\to G$ of two finitely generated subgroups $E', E''$ of $F_2$
onto a finitely presented group $G$, that is $E=\{(u,v)\in E'\times E'' \ | \ \phi(u)=\psi(v)\}$. Since every finitely generated subgroup of $F_2$ has bounded distortion, $E'\times E''$ is quasi-isometrically embedded into $F_2\times F_2$. So it is enough to show that
the distortion function of the equalizer of two homomorphisms $\phi: F_m\to G$ and
$\psi: F_n\to G$ of two free groups onto a finitely presented group $G$ in $F_m\times F_n$ is equivalent to the Dehn function of $G$.

Let $d(k)$ be the Dehn function of $G$. Let
$F_m=\langle x_1, ..., x_m\rangle$, $F_n=\langle y_1,...,y_n\rangle$ (we assume that these generating sets are closed under taking inversese). As a
generating set for $H=F_m\times F_n$ we take the set of all pairs $(x_i,1)$, $(1,y_j)$.
Let $E$ be the equalizer of $\phi$ and $\psi$ in $H$.

Let $r_1,...,r_\ell$ be generators of the kernel of $\psi$ (as a normal
subgroup of $F_n$). Without loss of generality we assume that the set $\{r_1,...,r_\ell\}$ is closed under cyclic shifts and inverses.

For every $i=1,...,m$ pick one word $t_i\in F_n$
such that $\phi(x_i)=\psi(t_i)$. For every $j=1,...,n$ pick one word
$s_j$ such that $\phi(s_j)=\psi(y_j)$. Then the equalizer $E$ is
generated by the pairs $(x_i,t_i), (s_j, y_j), (1,r_k)$, $i=1,...,m$,
$j=1,...,n$, $k=1,...,\ell$. Indeed, if $(u,v)\in E$, that is
 $\phi(u)=\psi(v)$
and $u=x_{i_1}x_{i_2}...x_{i_p}$ then
$$(u,v)=(x_{i_1},t_{i_1})...(x_{i_p}, t_{i_p})(1,a)$$
where $\psi(a)=1$. Since $a$ belongs to the kernel of $\psi$, we have
that $a$ is a product of conjugates of $r_k$:
$$
a=\prod_{i=1}^d r_{k_i}^{w_i}.
$$
Therefore
\begin{equation}
(u,v)=(x_{i_1},t_{i_1})...(x_{i_p}, t_{i_p})
\prod_{i=1}^d (1,r_{k_i})^{(w_i({\vec s}),w_i)}.
\label{11}
\end{equation}
Here $w({\vec s})$ denotes the word $w$ where each $y_j$ is substituted
by the corresponding $s_j$.

We shall prove that the distortion function of $E$ in $H$ is equivalent to $d(k)$.
In order to do that we need the following general statement.

\begin{lemma} Let $\Delta$ be a \vk diagram over a presentation $\langle X\ 
|\ R\rangle$ where $X=X\iv$, $R$ is closed under cyclic shifts and inverses. 
Let $w$ be the boundary label of $\Delta$.  Then $w$ is equal in the free 
group to a word of the form $u_1r_1u_2r_2...u_dr_du_{d+1}$ where:
\begin{enumerate} \item $r_i\in R$;
\item $u_1u_2...u_{d+1}=1$ in the free group; \item 
$\sum_{i=1}^{d+1}|u_i|\le 4e$ where $e$ is the number of edges of $\Delta$.
\end{enumerate} \label{lmm1}
\end{lemma}

{\bf Proof.} If $\Delta$ has an internal edge (i.e. an edge which belongs to 
the contours of two cells) then it has an internal edge $f$ one of whose 
vertices belongs to the boundary. Let us  cut $\Delta$ along $f$ leaving the 
second vertex of $f$ untouched. We can repeat this operation until we get a 
diagram $\Delta_1$ which does not have internal edges. It is easy to see 
that the boundary label of $\Delta_1$ is equal to $w$ in the free group. The 
number of edges of $\Delta_1$ which do not belong to contours of cells (let us call
them {\em edges of type 1} is the 
same as the number of such edges in $\Delta$ and the 
number of edges which belong 
to contours of cells in $\Delta_1$ ({\em edges of type 2})
is at most twice the number of such edges of 
$\Delta$ (we cut each edge from a contour of a cell 
at most once, after the cut we get two external 
edges instead of one internal edge).

Suppose that a cell $\Pi$ in $\Delta_1$ has more 
than one edge which has a 
common vertex with $\Pi$ but does not belong to the contour of $\Pi$ .  
Take any point $O$ on $\partial(\Pi)$. Let $p$ be the boundary path of $\Delta_1$
starting at $O$ and let $q$ be the boundary path of $\Pi$ starting at $O$.
Consider the path $qq\iv p$. The subpath $q\iv p$ bounds a subdiagram of 
$\Delta_1$ containing all cells but $\Pi$. Replace the path $q$ in $qq\iv p$
by a loop $q'$ with the same label starting at $O$ and lying inside the
cell $\Pi$. Let the region inside $q'$ be a new cell $\Pi'$. Then the path
$q'q\iv p$ bounds a diagram whose boundary label free is freely equal to $w$.
Notice that $\Pi'$ has exactly one edge having a common vertex with $\Pi'$
and not belonging to the contour of $\Pi$. Thus this operation reduces 
the number of cells which have more
than one edge which has a 
common vertex the cell but does not belong to the contour of it.  

After a number of such transsformations 
we shall have a diagram $\Delta_2$ 
which has the form of a tree $T$ with cells hanging like leaves (each has exactly 
one common vertex with the tree). 

The number of 
edges of type 1 in $\Delta_2$ cannot be bigger
than the number of all edges in $\Delta_1$, so it cannot be more 
than two times 
bigger than the total number of edges in $\Delta$. 

The boundary label of $\Delta_2$ is freely equal to $w$, and it has the form  
$u_1r_{i_1}u_2r_{i_2}...u_dr_{i_d}u_{d+1}$ where $d$ is the number of cells 
in $\Delta$, $u_1u_2...u_{d+1}$ is the label of a tree, so 
$u_1u_2...u_{d+1}=1$ in the free group. The sum of lengths 
of $u_i$ is 
at most four times the number of edges in $\Delta$ because
the word $u_1u_2...u_{d+1}$  is written on the tree $T$, and when we
travel along the tree,  we pass through each edge twice. 

The lemma is proved.

\medskip Let us consider the distortion function of $E$ in $H$.
Without loss of generality we can assume that $d(k)$ is the Dehn function of the presentation
$\langle y_1,...,y_n\ | \ r_1,...,r_\ell\rangle$ (recall that Dehn functions of different finite presentations of the same group are equivalent).

Let $(u,v)$ be any element in $E$ whose length in $H$, $|u|+|v|$, is $k\ge
 1$. Then as before
$$(u,v)=(x_{i_1},t_{i_1})...(x_{i_p}, t_{i_p})(1,a)$$
where $u=x_{i_1}...x_{i_p}$, $\psi(a)=1$.

Notice that the length of the word $a$
does not exceed $$|v|+c_1|u|\le c_1(|u|+|v|)=c_1k$$ where $c_1$
is the maximal length of $|t_i|$, $i=1,...,m$.

Let $\Delta$ be the minimal area \vk diagram over the presentation
$$\langle y_1,...,y_n\ |\ r_1,..., r_\ell\rangle$$ of $G$
with the boundary label $a$. Then the area of $\Delta$ does not exceed $d(c_1k)$. Since $\Delta$
is a planar graph, its number of edges $e$ does not exceed a constant times the area plus the length of the boundary of $\Delta$.

By Lemma \ref{lmm1}, $$a=u_1r_{i_1}u_2...r_{i_q}u_{q+1}$$ where $u_1...u_{q+1}=1$ in the free group, and $\sum |u_i|\le 4e$. Then
$$\begin{array}{l}
(1,a)=(u_1({\vec s})\cdot 1\cdot u_2({\vec s})\cdot 1\cdot ...\cdot 1\cdot u_{q+1}({\vec s}),\  u_1r_{i_1}u_2...r_{i_q}u_{q+1})\\
=u_1'\cdot (1,r_{i_1})\cdot u_2'(1,r_{i_2})\cdot...\cdot (1,r_{i_q})\cdot u_{q+1}'
\end{array}$$
where $u_i'$ denotes the word $u_i$ with letters $y_j$ substituted by $(s_j,y_j)$. Therefore
the length in $E$ of the element $(1,a)$ does not exceed $c_2d(c_1k)+c_1k$ for some constant $c_2$.
Hence the length in $E$ of the element $(u,v)$ does not exceed $k(1+c_1)+c_2d(c_1k)$.

This implies that the
distortion function of $E$ in $H$ does not exceed a function equivalent to
the Dehn function $d$.

To prove that the Dehn function $d$
does not exceed a function equivalent to the distortion
function of $E$ in $H$, it is enough, for every number $p\ge 1$, to take a word $a\in F_n$,
from the kernel of $\psi$, $|a|\le p$,  of area $d(p)$. Then it is easy to
see that
any representation of $(1,a)$ as a product of generators of
$E$ must contain at least $d(p)$ factors of the form $(1,r_k)$ (because it corresponds to a
representation of $a$ in the form $u_1r_{i_1}u_2r_{i_2}...r_{i_d}u_{d+1}$ where $u_1...u_{d+1}=1$ in the free group).
$\Box$

\subsection{Dehn Functions of Groups}

Our first goal is to give an almost complete description of
Dehn functions of finitely presented groups in terms of time functions of
Turing machines. First of all Birget and Sapir \cite{SBR}
proved that every Dehn function is the time
function of a nondeterministic Turing machine.

\begin{theorem} (Birget, Sapir, \cite{SBR})
Every Dehn function of a finitely presented group $G$ is equivalent to
the time function of some (not necessarily deterministic) Turing machine
solving the word problem in $G$.
\label{th111}
\end{theorem}

This result restricts the class of functions which can be Dehn functions of
groups. Indeed, time functions of non-deterministic machines are functions
$f(n)$ which can be computed deterministically in time at most $2^{f(n)}$.
It is easy to construct a recursive number $\alpha>2$ such that the function
$[n^\alpha]$ is not computable even in double exponential time by a Turing
machine, so $n^\alpha$ is not equivalent to the Dehn function of any
finitely presented group. This answers a question by Gersten (he asked if
every increasing recursive function $>n^2$ is equivalent to the Dehn
function of a finitely presented group).

\bigskip

The set of Dehn functions ``must" satisfy a yet another restriction: every
Dehn function ``must" be superadditive (more precisely, it ``must" be
equivalent to a superadditive function), that is $f(m+n)\ge f(m)+f(n)$ for
every $m,n$.  We put the word ``must" in quotation marks because the proof
of this restriction is yet to exist. Here is a quasi-proof. Notice that it
is enough to show that $f(m+n+c)\ge f(m)+f(n)$ for some constant $c$ and all
$m,n$ (since we identify equivalent functions).  Now, if the word $w$ of
length $\le m$ has area $f(m)$ and the word $w'$ of length $\le n$ has area
$f(n)$ and there are no cancellations in the product $w^gw'$ where $g$ is a
word of small length ($\le c/2$), then this product ``cannot" have area
smaller than $f(m)+f(m)$. Since the length of $w^gw'$ is $\le m+n+c$, we
have $f(m+n+c)\ge f(m)+f(n)$.  \setcounter{pdfive}{\value{ppp}} Figure
\thepdfive\  shows a diagram with boundary label $w^gw'$.

\bigskip
\unitlength=1mm
\linethickness{0.4pt}
\begin{picture}(126.67,38.00)
\put(46.83,20.50){\oval(47.00,29.67)[]}
\put(70.33,20.67){\line(1,0){9.00}}
\put(74.33, 23.67){\makebox(0,0)[cc]{$g$}}
\put(103.17,20.50){\oval(47.00,29.67)[]}
\put(46.00,38.00){\makebox(0,0)[cc]{$w$}}
\put(102.00,38.00){\makebox(0,0)[cc]{$w'$}}
\end{picture}
\begin{center}
\nopagebreak[4]
Fig. \theppp.

\end{center}
\addtocounter{ppp}{1}

Of course the problem is that we can probably tessellate the disk with the
boundary label $w^gw'$ in a different, more economical, way. Still there are
so many ways to choose $g$, and to connect two van Kampen diagrams that it
seems unlikely that we cannot find a product $w^gw'$ with area $f(m)+f(n)$.

Although, as we have said the proof of superadditivity property does not
exist at that time, Guba and Sapir were able to prove the following partial
result.

\begin{theorem} (Guba, Sapir, \cite{freeproducts}) The Dehn function of
every group which is a free product of two non-trivial groups, is
superadditive.  \label{freeproducts} \end{theorem}

The proof of this theorem basically shows that the idea presented above works
in the case of free products.

In view of Theorem \ref{freeproducts}, the superadditivity property is
equivalent to the following property:

\begin{quotation}
The Dehn function of any finitely presented group $G$ is equivalent to the
Dehn function of the free product $G * {\bf Z}$.
\end{quotation}

The next theorem gives a description of Dehn functions. It shows
that the class of functions $f(n) > n^4$ satisfying  restrictions
mentioned above virtually coincides with the class of Dehn functions $>n^4$ of
finitely presented groups.

\begin{theorem} (Sapir, Birget, Rips, \cite{SBR}) \label{1}
Let $M$ be a not necessarrily
deterministic Turing machine with time function
$T(n)$ for which $T(n)^4$ is superadditive.  Then there exists a finitely
presented group $G(M) = \langle A \rangle $ with Dehn function equivalent to
$T(n)^4$, and the smallest isodiametric function equivalent to $T(n)^3$.

Moreover, $G(M)$ simulates $M$, that is there exists an injective map $K$ from
the set of input words of $M$ to $(A \cup A^{-1})^+$ such that

\begin{enumerate}
\item $1/C|u|< |K(u)|<C|u|$ for some constant $C>1$ and
for every input word $u$;
\item An input word $u$ is accepted by $M$ if and only if $K(u)=1$ in $G$;
\end{enumerate}
\end{theorem}

This theorem implies the following description of the ``isoperimetric
spectrum" in $[4,\infty)$, that is the numbers $\alpha>4$ such that
$n^\alpha$ is equivalent to a Dehn function of a finitely presented group.

We say that a real number $\alpha$ is {\em computable in time} $\le T(m)$
for some function $T(m)$ if there exists a deterministic
Turing machine which for every
number $m$ written in binary computes the first $m$ digits of $\alpha$ in
time at most $T(m)$.

\begin{theorem} (Sapir, \cite{SBR})
\label{cy304}  For every real number $\alpha\ge 4$ computable in time
$\preceq 2^{2^m}$ the function $n^\alpha$ is equivalent to
the Dehn function of a finitely presented group and the smallest
isodiametric function of this group is $n^{3/4\alpha}$. On the other hand if
$n^\alpha$ is the Dehn function of a finitely presented group then
$\alpha$ is computable in time $\preceq 2^{2^{2^m}}$.
\end{theorem}

Of course all well known numbers $>4$ (say, rational numbers, $e+2, e\pi,
2\log_2\frac{a}{b}$ for integers $a$, $b$, $a>4b$), are computable in
polynomial time, so for these numbers $\alpha$, $n^\alpha$ is the Dehn
function of a finitely presented group. For $\alpha\le 4$, Brady and Bridson
proved that the spectrum contains all numbers of the form $2\log_2
\frac{2a}{b}$ where $a > b$ are integers, so the spectrum is
dense in the set of all real numbers, but
a description similar to Theorem \ref{cy304} is not known for numbers $\le
4$. Even for non-integer rational numbers between $2$ and $4$ we do not yet
know if they belong to the isoperimetric spectrum. We expect the result for
$\alpha \in [2, 4)$ to be similar to Theorem \ref{cy304}.

Of course Theorem \ref{1} provides examples of Dehn functions which are much
more complicated than $n^\alpha$. For example, functions like
$n^{2\pi}(\log n)^{\log n}\log\log n$ are clearly equivalent to fourth powers
of time functions of Turing machines (hint: take the Turing machines which
calculates the fourth root of such a function in the unary notation), and
by Theorem \ref{1} they are equivalent to Dehn functions of finitely
presented groups.

Theorem \ref{mikh} allows us to formulate the following corollary of Theorem
\ref{1} which gives examples of subgroups of the direct product of two free
groups with ``arbitrary weird" distortion.

\begin{theorem} (Sapir) For every time function $T(n)$ of a
non-deterministic Turing machine with superadditive $T(n)^4$ there exists a
subgroup of $F_2\times F_2$ with distortion function $T(n)^4$. In particular
for every real number $\alpha\ge 4$ computable in time $\le 2^{2^m}$ there
exists a subgroup of $F_2\times F_2$ with distortion function equivalent to
$n^\alpha$.  \label{sd}
\end{theorem}

Recall that $F_2 \times F_2$ is automatic. Notice that every cyclic
subgroup of it has (obviously) bounded distortion.

\subsection{Length Functions of a Finitely Generated Group}
\label{lf}

Theorems \ref{mikh} and \ref{sd} gives information about the set of distortion
functions of subgroups of one particular group, $F_2\times F_2$. In this
section, we shall fix an arbitrary finitely generated group $G$ and describe
all possible length functions (and hence distortion functions) of $G$ inside
other groups.

A complete description of all length functions of a finitely generated group
is given by the following theorem.

\begin{theorem} (Olshanskii, \cite{Ol4}) Let $\ell: G\to {\bf N}$
be a length function on a group $G$. Then the following conditions hold:
\begin{enumerate}
\item[(D1)] $\ell(g)=\ell(g\iv)$ for every $g\in G$; $\ell(g)=0$ if and only if $g=1$.
\item[(D2)] $\ell(gh)\le \ell(g)+\ell(h)$ for every $g,h\in G$.
\item[(D3)] There exists a positive number $c$ such that the cardinality of the set $\{g\in G\ | \ \ell(g)\le r\}$ does not exceed $c^r$ for every $r\in {\bf N}$.
\end{enumerate}
Conversely for every group $G$ and every function $\ell: G\to {\bf N}$ satisfying (D1) - (D3),
there exists an embedding of $G$ into a 2-generated group $H$ with  generating set $B=\{b_1,b_2\}$ such that the length function $g\to |g|_B$ is equivalent to $f$.
\label{ol4}
\end{theorem}

In the particular case when $G$ is a cyclic group,  Theorem \ref{ol4}
implies that for any number $\alpha \in (0, 1]$, there exists a group
$H_\alpha>G$ and an element $g\in H_\alpha$ such that the length of $g^i$
in $H_\alpha$ grows as $i^\alpha$.  This gives an answer to Gromov's
question \cite{gromov1}.

Another problem by Gromov \cite{gromov1} asked for a description of length
functions of cyclic groups in finitely presented groups. It is clear that
not every function satisfying (D1)--(D3) can be a length function of the
cyclic group in a finitely presented group: the cardinality of the set of
$O$-equivalence classes of functions satisfying (D1)-(D3) is continuum,
and the set of embeddings of
the infinite cyclic group into finitely presented groups is countable.
Nevertheless the following theorem shows that all
``reasonable" functions are length functions of a given finitely generated
group $G$ in a finitely presented group.

Let $G$ be a group with a finite generating set $A = \{a_1,\dots,a_m\}$.
Let $F_m$ be the free group generated by $A\cup A\iv$.  Every function
$\ell: G\to {\bf N}$ can be naturally extended to a function $\ell^*: F_m\to
N$.
We say that $\ell$ is {\em computable} if $\ell^*$
is computable in the natural sense.

\begin{theorem} (Olshanskii, \cite{Ol1}) Let $\ell$ be a computable function
$G\to {\bf N}$ satisfying (D1)-(D3). Then $G$ can be embedded into a
finitely presented group $H$ in such a way that the corresponding length
function is equivalent to $\ell$.  \label{ol}
\end{theorem}

This theorem immediately follows from Theorem \ref{ol4} and the following
result.

\begin{theorem} (Olshanskii,\cite{Ol1})
Every finitely generated and recursively presented group $G$ can be
quasi-isometrically embedded into a finitely presented group.  \label{ol1}
\end{theorem}

Although Theorem \ref{ol}  shows that all ``reasonable" functions are length
functions of a given finitely generated recursively presented group inside
finitely presented groups, it does not give a characterization of these
functions. Such a characterization has been found recently by Olshanskii.
This answers questions asked by P. Papasoglu and R. Gilman. It also gives a
complete solution of Gromov's problem from \cite{gromov1}.  It turned out
that such a characterization can be easily deduced from \cite{Ol4} and
 \cite{Ol1}.

We say that a function $\ell: G\to {\bf N}$ satisfies condition (D4) if
there exists a natural number $n$ and a recursively enumerable set $S\subset
F_m\times F_n$ such that

(a) if $(v_1,u), (v_2,u)\in S$ for some words $v_1,v_2,u$ then
$v_1$ and $v_2$ represent the same element in $G$;

(b) $\ell^*(v)=\min(\{|u|\; |\ (v,u)\in S\})$ for every $v\in F_m$.

Clearly it does not depend on the choice of generators of $G$ whether $\ell$
satisfies condition (D4) or not because of the obvious rewriting.

Notice that in (D4), we can always assume $n=2$. Indeed, if condition (D4)
holds for a function $\ell$ and some $n$, it also holds for $\ell$ and any
natural number $n'\ge 2$ since there is an isomorphic embedding of $F_n$
into $F_{n'}$.

\begin{theorem} (Olshanskii, \cite{os2}) Let $G$ be a finitely generated subgroup of a finitely
presented group $H$. Then the corresponding length function on $G$ satisfies
conditions (D1)--(D4). Conversely, for every finitely generated group $G$
and every  function $\ell: G\to {\bf N}$ satisfying conditions (D1)--(D4),
there exists an embedding of $G$ into a finitely presented group $H$ such
that the length function $g\mapsto |g|_H$ is $O$-equivalent to $\ell$.
\label{ol98} \end{theorem}

Condition (D4) is relatively complicated. We do not know if it is possible
to simplify it in general. But in the case when the group $G$ has solvable
word problem, including the important case when $G$ is cyclic, condition
(D4) can be replaced by a much simpler condition.

As usual, the {\em graph} of a function $\ell^*: F_m\to {\bf N}$ is the set
$(w, \ell^*(w))\subseteq F_m\times N$. A pair $(w,k)$ is said to lie {\em
above} the graph of $\ell^*$ if $\ell^*(w)\le k$.

\begin{theorem} (Sapir, \cite{os2}) Let $G$ be a finitely generated group
with decidable word problem. Then the function $\ell: g\mapsto |g|_H$ given
by an embedding of $G$ into a finitely presented group $H$ satisfies
condtions (D1)--(D3) and the following condition

(D4') The set of pairs above the graph of $\ell^*$ is recursively enumerable.

Conversely, for every function $\ell: G\to {\bf N}$ satisfying conditions
(D1), (D2), (D3), and (D4'), there exists an embedding of $G$ into a finitely
presented group $H$ such that the corresponding length function on $G$ is
$O$-equivalent to $\ell$.
\label{os98}
\end{theorem}

It is again clear that whether condition (D4') holds or not does not depend
of the choice of generators of $G$.

In the important particular case when $G$ is the infinite cyclic group we
have

     \begin{cor} (1)Let $g$ be an element of infinite order
     in a finitely presented group $H$ with a generating set
     ${\cal B} = \{b_1,\dots,b_k\}$. Denote
     $\ell(i)=|g^i|_{\cal B} =|g^i|$ for $i\in {\bf Z}$. Then
     \begin{itemize}

     \item[(C1)] $\ell(i)=\ell(-i)$ for $i\in{\bf Z}$ ($l$ is symmetric), and $\ell(i)=0$ iff $i=0$;

     \item[(C2)] $\ell(i+j)\le \ell(i)+\ell(j)$ for $i,j\in {\bf Z}$
     ($l$ is subadditive);

     \item[(C3)] there is a positive number $c$ such that
     $\,card\{i\in {\bf Z}| \ell(i)\le r\}\le c^r$ for any $r\in{\bf N}$.

     \item[(C4)] the set of integer pairs above the graph of $\ell$ is recursively enumerable.
\end{itemize}
     (2) Conversely, for any function $\ell:\,{\bf Z}\rightarrow{\bf N}$,
     satisfying the conditions (C1)--(C4), there is a
     finitely presented group $H$ and an element $g\in H$ such that
     $|g^i|_H$ is $O$-equivalent to $\ell(i)$.
\label{c}
\end{cor}

It is easy to prove that  (D4) implies (D4').  Indeed, suppose that (D4)
holds. Consider a Turing machine $M$ listing elements of the recursively
enumerable set $E$.  Let us change the machine $M$ in such a way that (1)
instead of pairs $(w,u)$ from $F_m\times F_n$ it produces pairs $(w, |u|)$
from $F_m\times {\bf N}$ and  (2) after every, say, 10, steps of
calculation, it goes through all pairs listed so far and for each of these
pairs $(w_i, k_i)$ adds a pair $(w_i, k_i+1)$ to the list, then it does the
next 10 steps of calculations, etc. Clearly, this new machine will list all
pairs which are above the graph of $\ell^*$ and only these pairs.
Thus the set of pairs above the
graph of $\ell^*$ is recursively enumerable and condition (D4') holds.

By the proper choice of a universal group $H$ it is not difficult to sharpen
the formulation of Theorems \ref{ol} and \ref{ol98}. One can select the
group $H$ in these theorems (independently of $G$) as the receptacle of all
possible ``computable distortions" of finitely generated recursively
presented groups. The next theorem follows from Theorem 4 from \cite{Ol1}.

\begin{theorem} (Olshanskii, 1998) There exists a finitely presented group $H$,
having the following property. For an arbitrary finitely generated
recursively presented group $G$ and an arbitrary function $\ell: R\to {\bf
N}$ satisfying conditions (D1)-(D4) there exists an embedding of $G$ into
$H$ such that the length function of $G$ corresponding to this embedding is
$O$-equivalent to $\ell$.
\end{theorem}

\subsection{Groups with Word Problem in NP}

The well known Higman theorem says that a group has a recursive presentation
if and only if this group is embeddable into a finitely presented group.
Theorem \ref{ol1} strengthens this result. The next Theorem strengthens it
even further.

\begin{theorem} (Birget, Olshanskii, Rips, Sapir \cite{BORS}) Let $G$ be a
finitely generated group with word problem solvable by a non-deterministic
Turing machine with time function $\le T(n)$ such that $T(n)^4$ is
superadditive. Then $G$ can be quasi-isometrically embedded into a finitely
presented group $H$ with isoperimetric function equivalent to $n^2T(n^2)^4$.
In particular, the word problem of a finitely generated group is in NP if
and only if  this group is a (quasi-isometric) subgroup of a finitely
presented group with polynomial isoperimetric function.  \label{th1}
\end{theorem}

\bigskip
In particular, this theorem gives a Higman-like  description of groups with
word problem in NP.

The class of finitely generated groups with word problem in NP is very
large. It clearly includes all matrix groups over ${\bf Q}$.
It also includes

\begin{itemize}
\item All finitely generated matrix groups over arbitrary fields: this
follows from the fact that every finitely generated field is a finite
extension of a purely transcendental extension of its simple subfield, and
the fact that the word problem in the ring of polynomials over ${\bf  Q}$ or
${\bf Z}/p{\bf Z}$ is solvable in polynomial time,
\item Polycyclic and finitely generated metabelian groups because they are
representable by matrices \cite{Segal},
\item Automatic groups (in particular, hyperbolic groups) \cite{WPG},
\item Groups of piecewise linear transformations of a line with finitely
many rational singularities (including the R. Thompson group $F$) \cite{BS},
\item Every finitely generated subgroup of a diagram group \cite{GS},
\item Every free Burnside group $B(m,n)$ for sufficiently large odd
exponent $n$ (see, for example, Storozhev's
argument in Section 28 of \cite{Ol89}).
\end{itemize}

This class is closed under free and direct products. It is easy to see using
Magnus' embedding that for every normal subgroup $N$ of a free finitely
generated group $F$ if $F/N$ has word problem in NP (resp. P) then $F/N'$
has word problem in NP (resp. P). Therefore every free group in the variety
of all solvable groups of a given class has word problem in P.

\bigskip

It is an interesting question whether this class also contains all
one-relator groups. There are of course finitely generated groups with word
problem not in NP, for example groups with undecidable word problem.
Moreover the construction from \cite{SBR} allows one to construct groups
with decidable but arbitrary hard word problem.  But these groups are in
some sense ``artificial". So perhaps the class of groups with word problem
in NP (which by Theorem \ref{th1} is the class of all subgroups of finitely
presented groups with polynomial Dehn functions) can be considered as the
class of ``tame" groups.

\bigskip

An example of an embedding of one group into another where lengths are not
distorted but areas are distorted can be found in Gersten \cite{Gersten3}.
Some examples of groups with big Dehn functions embeddable into groups with
small Dehn functions can be found in Madlener, Otto \cite{MO} and Baumslag,
Bridson, Miller and Short \cite{BBMS}.  Our results show that any
recursively presented finitely generated group can be embedded into a
finitely presented group with bounded length distortion but with close to
maximal possible area distortion.  Indeed, Theorem \ref{th111} shows that
an isoperimetric function of a group $H$ containing a given group $G$
cannot be smaller than the time complexity $T(n)$ of the word problem for
$G$, and Theorem \ref{th1} shows that $G$ can be embedded into a finitely
presented group with Dehn function at most $n^2T(n^2)^4$ (which is
polynomially equivalent to $T(n)$ ).

For matrix groups our theorem implies that every such group is embedded
quasi-isometrically into a finitely presented group with Dehn function
at most $n^{10+\epsilon}$ for every $\epsilon>0$.  It is interesting to know
the smallest Dehn function of a finitely generated group containing, for
example, the Baumslag-Solitar group $BS_{2,1}$.

\medskip

Notice that a semigroup analog of
Theorem \ref{th1} was obtained in Birget \cite{Bi}.

\medskip

As it usually happens, solution of one problem leads to solutions of other
problems.

In 1976, D. Collins asked \cite{Kour} if there exists a version of the Higman
embedding theorem which preserves the degree of unsolvability of the conjugacy
problem. The answer is ``yes" as the following theorem shows.

\begin{theorem} (Olshanskii, Sapir, 1998) The embedding described in the
proof of Theorem \ref{th1} preserves the degree of unsolvability of the
conjugacy problem. In particular, the conjugacy problem is decidable in $G$
if and only if it is decidable in $H$.
\label{os981}
\end{theorem}

Using the proof of Theorem \ref{th1}, in order to embed a finitely generated
group $G$ with word problem in NP into a finitely presented group with
polynomial isoperimetric function, one needs first construct a Turing
machine which solves the word problem in $G$, then convert it into a
so called $S$-machine (see below), then convert the $S$-machine into a
group. As a result the group we construct will have a relatively complicated
set of relations.  In some important cases like the Baumslag-Solitar group
$G_{2,1}$, the free Burnside groups $B(m,n)$, where $n$ is odd and $>>1$,
and others, we can modify our construction and get simple presentations of
groups with polynomial isoperimetric functions where these groups embed.

Consider, for example, the free Burnside group $B(m,n)$ with $m$ generators
$\{a_1,...,a_m\}$ and exponent $n$. This group is
very complicated and in particular not finitely presented if $m\ge 2$ and
$n$ is odd and $\ge 665$ (Adian, \cite{Adi}).  Now we are going to give a
presentation of a finitely presented group $H$ with a polynomial
isoperimetric function, containing $B(m,n)$ as a quasi-isometric subgroup.

The relations of $B(m,n)$ have the form $u^n=1$ where $u$ is an arbitrary
word in the alphabet of generators. So our goal is to find a finite set of
relations of a bigger group which will imply all the relations $u^n=1$ (and
no extra relations between generators of $B(m,n)$).

Instead of first writing relations of $H$, and then drawing van Kampen
diagrams we shall first draw diagrams, and then write relations.

For simplicity take $n=3$. The construction really does not depend much on
$n$, so we shall sometimes write $n$ instead of $3$. First of all, we shall
find a finite set of relations which imply relations of the form
$$K(uq_1uq_2uq_3)=k_1(uq_1uq_2uq_3)k_2(uq_1uq_2uq_3)'k_3....k_N(uq_1uq_2uq_3)^{(N)}$$ for every word $u$ in the alphabet
$\{a_1,...,a_m\}$. Here $N$ is a sufficiently large number (28 is enough),
$k_1,...,k_N, q_1,q_2, q_3$ are new letters, and
the  words between consecutive
$k's$ are copies of $uq_1uq_2uq_3$ written in disjoint alphabets.  The group
given by these relations will be denoted by $G_{m,n}$.
\setcounter{pdsix}{\value{ppp}} Figure \thepdsix\  shows the van Kampen
diagram (below it will be called a {\em disc}) with boundary label
$K(uq_1uq_2uq_3)$.

\bigskip

\unitlength=0.8mm
\linethickness{0.4pt}
\begin{picture}(140.00,92.00)
\put(89.17,46.00){\oval(101.00,83.33)[]}
\put(88.33,45.50){\oval(86.00,69.67)[]}
\put(88.67,44.50){\oval(73.33,56.33)[]}
\put(88.83,44.17){\oval(59.67,44.33)[]}
\put(89.33,43.50){\oval(48.67,32.33)[]}
\put(89.33,43.00){\oval(37.33,22.67)[]}
\put(90.67,42.83){\oval(25.33,13.00)[]}
\put(78.00,43.00){\line(-1,1){39.33}}
\put(44.00,87.67){\line(1,-1){38.33}}
\put(97.00,49.33){\line(1,1){38.00}}
\put(103.00,44.67){\line(1,1){36.67}}
\put(38.67,10.00){\line(4,3){40.33}}
\put(44.33,4.33){\line(6,5){39.00}}
\put(133.67,4.33){\line(-1,1){33.67}}
\put(100.00,38.00){\line(0,0){0.00}}
\put(139.67,10.00){\line(-6,5){37.00}}
\put(69.00,72.67){\vector(1,0){5.33}}
\put(97.33,72.67){\vector(1,0){6.00}}
\put(114.33,72.67){\vector(1,0){6.00}}
\put(68.00,72.67){\circle*{0.67}}
\put(96.33,72.67){\circle*{0.67}}
\put(113.67,72.67){\circle*{0.67}}
\put(69.00,80.33){\vector(1,0){5.33}}
\put(97.33,80.33){\vector(1,0){6.00}}
\put(68.00,80.33){\circle*{0.67}}
\put(96.33,80.33){\circle*{0.67}}
\put(120.67,80.33){\circle*{0.67}}
\put(68.00,72.67){\line(-5,6){6.33}}
\put(74.33,72.67){\line(0,1){7.67}}
\put(96.33,72.67){\line(-4,5){6.00}}
\put(103.33,72.67){\line(0,1){7.67}}
\put(39.33,87.33){\makebox(0,0)[cc]{$k_1$}}
\put(140.00,86.67){\makebox(0,0)[cc]{$k_2$}}
\put(139.67,3.33){\makebox(0,0)[cc]{$k_3$}}
\put(71.00,70.33){\makebox(0,0)[cc]{$q_1$}}
\put(99.00,70.67){\makebox(0,0)[cc]{$q_2$}}
\put(67.33,82.00){\makebox(0,0)[cc]{$aq_1$}}
\put(97.00,82.67){\makebox(0,0)[cc]{$aq_2$}}
\put(117.00,70.67){\makebox(0,0)[cc]{$q_3$}}
\put(120.67,82.67){\makebox(0,0)[cc]{$aq_3$}}
\put(82.33,92.00){\makebox(0,0)[cc]{$uq_1uq_2uq_3$}}
\put(113.67,72.67){\line(0,1){7.67}}
\put(120.67,80.33){\vector(1,0){3.67}}
\end{picture}
\begin{center}
\nopagebreak[4]
Fig. \theppp.

\end{center}
\addtocounter{ppp}{1}

On the boundary of this diagram we have the word $K(uq_1uq_2uq_3)$. The
words on each of the concentric circles is labeled by $K(u_iq_1u_iq_2u_iq_3)$
where $u_i$ is a prefix of $u$ of length $i-1$. The word written on the
innermost circle is
$K(q_1q_2q_3)$. This  word will be called the {\em hub}.
The edges connecting the circles are labeled by letters
$r_1,...,r_m$. The cells tessellating the space between the circles
have labels

\begin{itemize}
\item $q_i^{r_j}=a_jq_i$, $i=1,2,3$, $j=1,...,m$.
\item $ar=ra$, $a\in \{a_1,...,a_m\}$, $r\in  \{r_1,...,r_m\}$
\item $kr=rk$, $k\in\{k_1,...,k_N\}$, $r\in \{r_1,...,r_m\}$.
\end{itemize}

\noindent plus $N$ copies of each of these relations written in $N$ disjoint
alphabets.  These relations plus the {\em hub relation} $K(q_1q_2q_3)$ form
the presentation of $G_{m,n}$.

Now we construct $H=H_{m,n}$. Take a copy of $B(m,n)$ generated by
$\{b_1,...,b_m\}$.  The group $H$ will be an HNN-extension of the direct
product  $B(m,n)\times G_{m,n}$.  Here is the \vk diagram:

\vskip -0.6 in

\unitlength=0.8 mm
\linethickness{0.4pt}
\begin{picture}(134.67,113.00)
\put(79.50,28.83){\oval(63.00,46.33)[]}
\put(89.50,50.17){\oval(90.33,97.00)[lb]}
\put(98.83,50.83){\oval(32.33,98.33)[rb]}
\put(53.67,54.00){\oval(18.67,8.00)[lt]}
\put(108.17,51.50){\oval(13.67,13.00)[rt]}
\put(44.33,47.00){\line(0,1){7.33}}
\put(53.33,5.67){\line(-2,-3){2.67}}
\put(48.00,11.00){\line(-6,-5){3.67}}
\put(48.00,46.67){\line(-1,2){3.67}}
\put(53.33,52.00){\line(-5,6){5.00}}
\put(105.33,52.00){\line(3,5){3.67}}
\put(110.33,48.33){\line(4,3){4.67}}
\put(105.33,5.67){\line(1,-1){4.00}}
\put(111.00,11.67){\line(1,-1){4.00}}
\put(80.67,1.67){\line(1,0){23.67}}
\bezier{492}(53.33,58.00)(82.00,113.00)(108.33,58.00)
\bezier{180}(108.67,58.00)(94.67,78.33)(80.00,64.00)
\bezier{136}(108.67,57.33)(91.00,52.67)(80.00,64.00)
\put(78.67,89.00){\makebox(0,0)[cc]{$uq_1uq_2uq_3$}}
\put(90.67,65.00){\makebox(0,0)[cc]{$u_b^3$}}
\put(44.33,59.00){\makebox(0,0)[cc]{$k_1$}}
\put(114.33,57.67){\makebox(0,0)[cc]{$k_2$}}
\put(115.33,1.67){\makebox(0,0)[cc]{$k_3$}}
\put(42.33,2.33){\makebox(0,0)[cc]{$k_4$}}
\put(76.00,48.33){\makebox(0,0)[cc]{$uq_1uq_2uq_3$}}
\put(54.33,54.00){\makebox(0,0)[cc]{$\rho$}}
\put(76.67,25.67){\makebox(0,0)[cc]{Disc}}
\end{picture}
\setcounter{pdseven}{\value{ppp}}
\begin{center}
\nopagebreak[4]
Fig. \theppp.
\end{center}
\addtocounter{ppp}{1}

This is an annular diagram. The hole of it has label $u_b^n$ ($u_b$ is the
word $u$ rewritten in the alphabet $\{b_1,...,b_m\}$).  The boundary label
of the disc is $K(uq_1uq_2uq_3)$, the label of the external boundary of the
diagram is also $K(uq_1uq_2uq_3)$. In order to fill this diagram as shown on
the picture, one needs a new (stable) letter $\rho$ and the following
relations:

\begin{itemize}
\item $\rho k= k\rho$ for  $k \in \{k_1,...,k_N\}$.
\item $\rho q= q\rho$ for every $q\in \{q_1,q_2,q_3,...,q_3^{(N)}\}$.
\item $\rho a=a\rho$ for every $a$ from the $N$ copies of $\{a_1,...,a_m\}$.
\item $a^\rho=ab$ for every $a\in \{a_1,...,a_m\}$.
\item $ab=ba$ for every $a\in \{a_1,...,a_m\}$ and $b\in \{b_1,...,b_m\}$.
\item $q_ib=bq_i$ for every $b\in \{b_1,...,b_m\}$, $i=1,2,3$.
\end{itemize}

Since the label of the external boundary is $K(uq_1uq_2uq_3)$, we can glue
in a disc with this label, and make our annular diagram into an ordinary
diagram with boundary label $u_b^n$. Since the discs are filled with cells
corresponding to the relations of $G_{m,n}$, and the rest is filled with
cells corresponding to the new relations, we get that all defining relations
of $B(m,n)$ follow from the (finitely many) relations that we got. The
group $H$ that we just created is what we need.

\begin{theorem} (Olshanskii, Sapir, 1998) The natural homomorphism of $B(m,n)$ into $H$ is a
quasi-isometric embedding.
The group $H$ has isoperimetric function $n^{8+\epsilon}$  provided $n$ is odd and $\ge 10^{10}$; $\lim_{n\to \infty} \epsilon=0$.
\label{sss}
\end{theorem}

Similarly we can quasi-isometrically embed a relatively free group $G$ of
any finitely based group variety into a finitely presented group. The
resulting group will have a polynomial isoperimetric function provided $G$
has polynomial {\em verbal isoperimetric function}. This function is defined
as follows:

\begin{quotation}
Let $v(x_1,...,x_n)$ be a word. Suppose that $w$ is in the verbal subgroup
vsg$(v)$. Then $w=\prod_i v(X_{1,i},...,X_{n,i})$. Fix such a representation of $w$ with minimal sum of lengths of all $|X_{j,i}|$ involved in this representation. The verbal isoperimetric
function gives an upper bound for this sum in terms of $|w|$.
\end{quotation}

This function does not depend (up to ``big O")
on the identity defining the variety, so one
can speak about verbal isoperimetric functions of varieties.

For example, the variety of solvable groups has polynomial verbal
isoperimetric function, so our construction embeds it into a group with
polynomial Dehn function.

The variety of Burnside groups of odd exponent $n>>1$ has verbal
isoperimetric function $n^{1+\epsilon}$ ($\lim_{n\to \infty} \epsilon=0$).
This can be proven by modifying
Storozhev's argument from \cite{Ol89} (Storozhev's argument gives estimate
$n^4$ for the verbal isoperimetric function).

We can also embed in a similar way the
Baumslag-Solitar groups $G_{1,n}$ into finitely presented groups with
isoperimetric function $n^{10}$.

\section{Methods}

\subsection{$S$-machines}

First of all let us present some ideas how to find a group with an
``arbitrary" Dehn function.  Consider again the main diagram called a {\em
disc} on the Figure \thepdsix\ for the group $G_{m,n}$.

The disc is divided by the $k$-bands into $N$ sectors.  The words written on
the circles between consecutive $k$'s have the form

$$uq_1uq_2uq_3$$
and to pass from one level to another level we replace $q_i$ by $aq_i$.
So we can imagine these words written on a tape of a Turing machine, $q_i$ mark
the places where the heads are, and we have a rule of the form
$$[q_1\to aq_1, q_2\to aq_2, q_3\to aq_3]$$
for every $a$.

\bigskip

What we get is a simple example of a so called {\em $S$-machine}.

Roughly speaking, the difference between $S$-machines and ordinary Turing
machines is that $S$-machines are almost ``blind". They ``see" letters
written on the tape only when these letters are between two heads of the
machine and the heads are very close to each other.  If the heads are far
apart, the machine does not see any letters on the tape, in this case a
command executed by the machine depends only on the state of the heads.

In contrast, ordinary Turing machines can see letters on the tape near the
position where the head is. The command executed by the machine always
depends not only on the state of the head but (which is very important!)
also on the letter(s) observed by the head.  Notice that even for moving the
head a Turing machine one square to the left, one needs to know the content
of the square to the left of the head.

Let us give a precise definition of $S$-machines. Let $k$ be a natural
number.  Consider now a {\em language of admissible words}. It consists of
words of the form $$q_1u_1q_2...u_kq_{k+1}$$ where $q_i$ are letters from
disjoint sets $Q_i$, $i=1,...,k+1$, $u_i$ are reduced group words in an
alphabet $Y_i$ ($Y_i$ are not necessarily disjoint), the sets $\bar
Y=\bigcup Y_i$ and $\bar Q=\bigcup Q_i$ are disjoint.

Notice that in every admissible word, there is exactly one representative of
each $Q_i$ and these representatives appear in this word in the order of the
indices of $Q_i$.

If $0\le i\le j\le k$ and $W=
q_1u_1q_2...u_kq_{k+1}$ is an admissible word then
the subword $q_iu_i...q_j$ of $W$ is called the $(Q_i,Q_j)$-subword of $W$
($i<j$).

An $S$-machine is a rewriting system \cite{KharSap}. The objects
of this rewriting system are all admissible words.

The rewriting rules, or {\em $S$-rules}, have the following form:
$$[U_1\to V_1,...,U_m\to V_m]$$
where the following conditions hold:
\begin{description}
\item Each $U_i$ is a subword of an admissible word starting with
a $Q_\ell$-letter and ending with a $Q_r$-letter (where $\ell=\ell(i)$
must not exceed $r=r(i)$, of
course).
\item If $i<j$ then $r(i) < \ell(j)$.
\item Each $V_i$ is also a subword of an admissible word whose $Q$-letters
belong to $Q_{\ell(i)}\cup...\cup Q_{r(i)}$ and which contains a $Q_\ell$-letter
and a $Q_r$-letter.
\item If $\ell(1)=1$ then $V_1$ must start with a $Q_1$-letter and if
$r(m)=k+1$ then $V_n$ must end with a $Q_{k+1}$-letter (so tape letters
are not inserted to the left of $Q_1$-letters and to the right of $Q_{k+1}$-letters).
\end{description}

To apply an $S$-rule to a word $W$ means to replace simultaneously subwords
$U_i$ by subwords $V_i$, $i=1,...,m$. In particular, this means that our
rule is not applicable if one of the $U_i$'s is not a subword of $W$. The
following convention is important:

After every application of a rewriting rule, the word is automatically
reduced. We do not consider reducing of an admissible word a separate step
of an $S$-machine.

We also always assume that an $S$-machine is symmetric, that is for every
rule of the $S$-machine the inverse rule (defined in the natural way) is
also a rule of this $S$-machine.  This reflects the fact that the $r$-edges
in the disc on Figure \thepdsix\ can point away from the hub or toward the
hub.

Notice that virtually any $S$-machine is highly nondeterministic.

Among all admissible words of an $S$-machine we fix one word $W_0$. If an
$S$-machine $\sss$ can take an admissible word $W$ to $W_0$ then we say
that $\sss$ {\em accepts} $W$. We can define a time and space function of an
$S$-machine as usual. If $U\to U_1\to...\to U_n=W_0$ is an accepting
computation of the $\sss$-machine $\sss$ then $|U|+|U_1|+...+|U_n|$ is
called the {\em area} of this computation. This allows us to define the {\em
area function} of an $S$-machine.

\begin{theorem} (Sapir, \cite{SBR})  $S$-machines are polynomially
equivalent to Turing machines. More precisely for every Turing machine $M$
with time function $T(n)$ there exists an $S$-machine with area function
$T(n)^4$ which is equivalent to $M$ (this means that there exists a
correspondence $\phi$ between configurations of $M$ and admissible words of
$\sss$, given a configuration $c$, the word  $\phi(c)$ is computable in
linear time, and the machine $M$ accepts $c$ if and only if $\sss$ accepts
$\phi(c)$).  \label{sap}
\end{theorem}

In fact a stronger theorem can be deduced from the main results of
\cite{SBR}. It was recently proved by Sapir.

\begin{theorem} (Sapir, 1998) For every Turing machine $M$ with time
function $\le T(n)$ such that $T(n)^4$ is superadditive, there exists an
$S$-machine $\sss$ with one head and only one internal state which is
equivalent to $M$ and has time function $\le T(n)^4$. \label{sap98}
\end{theorem}

Notice that an $S$-machine with one head and one state letter is completely
blind (in the sense explained above). The rules of such an $S$-machine have
the following very simple form:

$$[q\to uqv]$$
where $q$ is the internal state, $u$ and $v$ are words in the tape alphabet.

The amazing fact is that the proof of a completely Computer Science
statement, Theorem \ref{sap98}, involves some heavy geometric group theory.
We first convert $M$ into an $S$-machine $\sss_1$ with many heads, then
convert  $\sss_1$ into the group from \cite{SBR} with Dehn function
$T(n)^4$, then convert the group into an $S$-machine $\sss$ with one head
and one internal state, having time function $T(n)^4$ (in the last step we
use an idea from Miller \cite{Miller}).

The group $G_N(\sss)$ associated with an $S$-machine $\sss$ is constructed
in the same way as the group $G_{m,n}$ presented above. We add all rules of
$\sss$ to the set of generators and for every rule $r$ of the form $[U_1\to
V_1,...,U_p\to V_p]$ we have $p$ relations $U_1^r=V_1,...,U_p^r=V_p$. These
relations replace the relations $q_i^r=aq_i$ in the presentation of
$G_{m,n}$. Other relations are the same.

Although this construction slightly differs from the construction in
\cite{SBR} it is possible to prove the following statement.

\begin{theorem}(Sapir, \cite{SBR}) Every Turing machine $M$ with time function $T(n)$ can be
converted into an $S$-machine $\sss$ in such a way that the Dehn function of
the group $G_N(\sss)$ is $T(n)^4$ provided $T(n)^4$ is superadditive.
\end{theorem}

Now in order to embed a finitely generated group $G$ into a finitely presented group
we take a Turing machine $M$ recognizing words which are equal to 1 in $G$, convert it into an $S$-machine $\sss$, and then basically repeat the construction of the group $H_{m,n}$
replacing $B(m,n)$ by $G$ and $G_{m,n}$ by $G_N(\sss)$. The resulting group is denoted by $H_N(\sss)$. It plays the role of group $H$ in Theorem \ref{th1}.

\subsection{Why $S$-machines?}

Here we will explain why we need to convert Turing machines into $S$-machines.

Consider any Turing machine $M$. For simplicity assume that $M$ has one tape, which is always finite, but
we can add squares at the right end of the tape, the alphabet $A$ of tape letters, the set $Q$ of states, and the set $R$ of transitions. As usual we assume that
the head is always placed between two squares of the tape, and observes both squares.
So the transitions have the form $uqv\to u'q'v'$ where $u, v, u', v'$ are words in the tape alphabet, $q,q'\in Q$ (see \cite{Ro} for details). Then using the same idea as in the construction of $G_{m,n}$ we can replace the relations $q_i^r=aq_i$ by
$(uqv)^r=u'q'v'$ (here $r$ is a letter associated with the transition of $M$). As in \cite{Ro} we assume that $M$ has only one accept configuration $W_0$.

Thus we have the following presentation of the group $G(M)$ associated with $M$.

\begin{itemize}
\item $(uqv)^r=u'q'v'$, for every transition $r=[uqv\to u'q'v']$ of the machine $M$,
\item $ar=ra$, for every $a\in A$ and every transition $r$
\item $kr=rk$, for every $k\in\{k_1,...,k_N\}$ and every transition $r$.
\end{itemize}

As before, we need $N$ copies of each of these relations written in disjoint alphabets. The hub relation will have the form $K(W_0)$ where $W_0$ is the accept configuration.

Now it is easy to see that for every accepted word $U$ we can tessellate the disc with boundary
label $K(U)$ into cells labeled by these relations. Let $U=U_1\to U_2\to ...\to U_p=W_0$ be an accepting computation.  As  before we will have a sequence of
concentric circles, each labeled by $K(U_i)$, the
innermost oval will be labeled by $K(W_0)$.

So it is easy to see that the word $K(U)$ is equal to 1 in $G(M)$ if the configuration $U$ is accepted by $M$.

Unfortunately the converse statement is wrong in most cases and this is precisely why we need $S$-machines. Let us demonstrate this on a simple example. Consider the following Turing machine $M$. It has two states $q, q_0$ and one tape letter $a$. The only transitions are the following:

\begin{itemize}
\item[($r_1$)] $aq\to q_0$,
\item[($r_2$)] $aq_0\to q_0$.
\end{itemize}

\noindent The stop configuration is $q_0$ (the tape is empty). It is clear that the set of configurations accepted by this machine consists of configurations $a^nq_0$ and $a^mq$ where
$n\ge 0, m>0$, and does not include, for example, the configuration $q$. Thus we would like $K(q)$ to be not equal to 1 in $G(M)$.
The
diagram on
\setcounter{pdeight}{\value{ppp}}
Figure \thepdeight\
shows that $K(q)$ is equal to 1 in this group.

\medskip
\unitlength=1.1mm
\linethickness{0.4pt}
\begin{picture}(109.00,66.67)
\put(76.17,33.83){\oval(65.67,60.33)[]}
\put(75.50,34.50){\oval(45.00,38.33)[]}
\put(75.83,33.67){\oval(20.33,16.00)[]}
\put(43.33,58.67){\line(1,-1){22.33}}
\put(48.67,64.00){\line(1,-1){22.33}}
\put(80.67,41.67){\line(1,1){22.00}}
\put(86.00,36.00){\line(1,1){22.67}}
\put(86.00,31.00){\line(6,-5){23.00}}
\put(80.33,25.67){\line(1,-1){22.67}}
\put(71.33,25.67){\line(-1,-1){22.33}}
\put(65.67,31.33){\line(-1,-1){22.33}}
\put(50.67,62.00){\vector(1,-1){3.67}}
\put(86.33,47.33){\vector(1,1){4.00}}
\put(98.33,59.33){\vector(-1,-1){3.67}}
\put(103.00,53.00){\vector(-1,-1){3.00}}
\put(91.33,41.33){\vector(1,1){4.00}}
\put(90.33,27.33){\vector(4,-3){4.33}}
\put(105.67,14.67){\vector(-4,3){4.33}}
\put(97.00,9.00){\vector(-1,1){3.33}}
\put(84.67,21.33){\vector(1,-1){4.33}}
\put(66.67,21.00){\vector(-1,-1){4.33}}
\put(60.00,25.67){\vector(-1,-1){3.67}}
\put(47.33,13.00){\vector(1,1){3.67}}
\put(54.00,8.33){\vector(1,1){4.00}}
\put(46.83,55.17){\vector(1,-1){3.67}}
\put(60.67,41.33){\vector(-1,1){3.67}}
\put(76.00,48.00){\line(-3,-4){4.67}}
\put(76.00,47.67){\circle*{1.33}}
\put(76.00,47.67){\line(0,1){6.00}}
\put(76.00,49.67){\vector(0,1){3.00}}
\put(76.00,47.84){\vector(-3,-4){3.33}}
\put(76.00,53.67){\line(0,1){5.00}}
\put(76.00,58.67){\line(-5,1){27.33}}
\put(73.33,53.67){\vector(-1,0){10.33}}
\put(74.00,41.67){\vector(1,0){5.67}}
\put(80.00,53.67){\vector(1,0){10.33}}
\put(72.00,64.00){\vector(1,0){18.67}}
\put(52.33,57.33){\makebox(0,0)[cc]{$r_1$}}
\put(63.67,46.00){\makebox(0,0)[cc]{$r_2$}}
\put(75.67,66.67){\makebox(0,0)[cc]{$q$}}
\put(83.67,55.67){\makebox(0,0)[cc]{$q_0$}}
\put(67.33,55.67){\makebox(0,0)[cc]{$a$}}
\put(75.33,39.00){\makebox(0,0)[cc]{$q_0$}}
\put(77.67,50.33){\makebox(0,0)[cc]{$r_2$}}
\put(73.00,46.67){\makebox(0,0)[cc]{$a$}}
\put(87.33,45.67){\makebox(0,0)[cc]{$r_2$}}
\put(70.00,61.67){\makebox(0,0)[cc]{$a$}}
\put(73.67,56.00){\makebox(0,0)[cc]{$r_1$}}
\put(76.00,58.33){\circle*{1.33}}
\put(76.00,57.67){\vector(0,-1){3.00}}
\put(63.83,61.00){\vector(-4,1){3.50}}
\put(68.83,43.83){\vector(-1,1){3.50}}
\put(67.00,40.33){\makebox(0,0)[cc]{$k_1$}}
\put(54.83,52.67){\makebox(0,0)[cc]{$k_1$}}
\put(44.50,62.67){\makebox(0,0)[cc]{$k_1$}}
\put(108.00,62.50){\makebox(0,0)[cc]{$k_2$}}
\put(96.83,52.17){\makebox(0,0)[cc]{$k_2$}}
\put(84.50,40.33){\makebox(0,0)[cc]{$k_2$}}
\end{picture}

\begin{center}
\nopagebreak[4]
Fig. \theppp.

\end{center}
\addtocounter{ppp}{1}

This picture shows the tessellation of only one of the sectors. The other sectors are tessellated in the same way.

One can easily see the difference between this picture and the standard picture of a disc. Here
we have pairs of cells which have two common edges, and in the standard disc cells could have at most one common edge.

The  diagram on
\setcounter{pdnine}{\value{ppp}}
Figure \thepdnine\
is a subdiagram of the diagram on Fig. \thepdeight. It
consists of two cells corresponding to the relations $r_2a=ar_2$ and $(aq_0)^{r_2}=q_0$ and has boundary label corresponding to the relation $(q_0)^{r_2}=a\iv q_0$ which is the
relation corresponding to the rule $[q_0\to a\iv q_0]$, the inverse rule for $[q_0\to aq_0]$.

\bigskip
\unitlength=1.00mm
\linethickness{0.4pt}
\begin{picture}(92.67,18.67)
\put(63.67,11.33){\vector(-1,1){3.33}}
\put(86.00,10.33){\vector(1,1){4.00}}
\put(75.67,11.00){\line(-3,-4){4.67}}
\put(75.67,10.67){\circle*{0.67}}
\put(75.67,10.67){\line(0,1){6.00}}
\put(75.67,12.67){\vector(0,1){3.00}}
\put(75.67,10.67){\vector(-3,-4){3.33}}
\put(73.00,16.67){\vector(-1,0){10.33}}
\put(73.67,4.67){\vector(1,0){5.67}}
\put(79.67,16.67){\vector(1,0){10.33}}
\put(63.34,9.00){\makebox(0,0)[cc]{$r_2$}}
\put(83.34,18.67){\makebox(0,0)[cc]{$q_0$}}
\put(67.00,18.67){\makebox(0,0)[cc]{$a$}}
\put(77.34,13.33){\makebox(0,0)[cc]{$r_2$}}
\put(72.67,9.67){\makebox(0,0)[cc]{$a$}}
\put(87.00,8.67){\makebox(0,0)[cc]{$r_2$}}
\put(58.67,16.67){\line(1,-1){12.00}}
\put(88.33,12.67){\line(1,1){4.33}}
\put(92.67,17.00){\line(-1,-1){12.33}}
\put(58.67,16.67){\line(1,0){34.00}}
\put(70.67,4.67){\line(1,0){9.67}}
\put(75.33,2.33){\makebox(0,0)[cc]{$q_0$}}
\end{picture}

\begin{center}
\nopagebreak[4]
Fig. \theppp.

\end{center}
\addtocounter{ppp}{1}

It is possible to prove that the group $G(M)$ actually simulates the $S$-machine with the set of admissible words $a^nq$, $a^nq_0$ and the set of rules

\begin{itemize}
\item $q\to a\iv q_0$,
\item $q_0\to a\iv q_0$.
\end{itemize}
plus the inverse rules.
This $S$-machine is ``stronger" than $M$, it accepts more configurations, including the configuration $q$.

In general if we take any Turing machine $M$ and repeat this construction we will get a group simulating the $S$-machine obtained by replacing every transition $uqv\to u'q'v'$ by the
$S$-rule $[q\to u\iv u'q'v'v\iv]$. This $S$-machine will almost always be much stronger than the original Turing machine.

One way around this problem was invented by Boone and Novikov \cite{Ro}. This is why they used the Baumslag-Solitar type relations $x^a=x^2$. These relations prevent appearance of negative letters on the ``tape" (the concentric circles in the disc). But we could not use these relations because they make the Dehn function exponential.

Thus we had to prove instead that $S$-machines are polynomially equivalent to ordinary Turing machines (Theorem \ref{sap}).

\subsection{Geometry of van Kampen Diagrams}

In order to analyze an arbitrary diagram over $H=H_{m,n}$, and $H_N(\sss)$
in general we change the presentation of
$H$. We add all words $K(u)$ (discs) and all relations of $G=\langle B\rangle$ to the presentation.
The presentation becomes infinite. After that we order the relations, saying
that the discs have the highest rank, $r$-relations have smaller ranks, and the $b$-commutativity relations have the lowest rank. With every diagram we associate its type, a vector,  the first coordinate of which is the number of discs, and the last coordinate
is the number of $b$-commutativity cells (we omit the ranking of other relations). It turns out that diagrams of minimal type have nice geometric properties.

The main and easy concept which helps us analyze these diagrams is the
concept of a {\em band} \footnote{Other people call them corridors and strips.}.
 If $S$ is a set of letters then an $S$-band is a
sequence of cells $\pi_1,...,\pi_n$ in a van Kampen diagram such that each
two consecutive cells in this sequence have a common edge labeled by a
letter from $S$. \setcounter{pdten}{\value{ppp}} Figure \thepdten\
illustrates this concept.

 \bigskip
\unitlength=0.90mm
\linethickness{0.4pt}
\begin{picture}(145.33,30.11)
\put(19.33,30.11){\line(1,0){67.00}}
\put(106.33,30.11){\line(1,0){36.00}}
\put(142.33,13.11){\line(-1,0){35.67}}
\put(86.33,13.11){\line(-1,0){67.00}}
\put(33.00,21.11){\line(1,0){50.00}}
\put(110.00,20.78){\line(1,0){19.33}}
\put(30.00,8.78){\vector(1,1){10.33}}
\put(52.33,8.44){\vector(0,1){10.00}}
\put(76.33,8.11){\vector(-1,1){10.33}}
\put(105.66,7.78){\vector(1,2){5.33}}
\put(132.66,8.11){\vector(-1,1){10.00}}
\put(19.33,30.11){\vector(0,-1){17.00}}
\put(46.66,29.78){\vector(0,-1){16.67}}
\put(73.66,30.11){\vector(0,-1){17.00}}
\put(142.33,30.11){\vector(0,-1){17.00}}
\put(117.00,30.11){\vector(0,-1){17.00}}
\put(16.00,21.11){\makebox(0,0)[cc]{$e$}}
\put(29.66,25.44){\makebox(0,0)[cc]{$\pi_1$}}
\put(60.33,25.44){\makebox(0,0)[cc]{$\pi_2$}}
\put(133.00,25.78){\makebox(0,0)[cc]{$\pi_n$}}
\put(145.33,21.44){\makebox(0,0)[cc]{$f$}}
\put(77.00,25.11){\makebox(0,0)[cc]{$e_2$}}
\put(122.00,25.44){\makebox(0,0)[cc]{$e_{n-1}$}}
\put(100.33,25.44){\makebox(0,0)[cc]{$S$}}
\put(96.33,30.11){\makebox(0,0)[cc]{$\dots$}}
\put(96.33,13.11){\makebox(0,0)[cc]{$\dots$}}
\put(26.66,4.78){\makebox(0,0)[cc]{$l(\pi_1,e_1)$}}
\put(52.66,4.11){\makebox(0,0)[cc]{$l(\pi_2,e_1)$}}
\put(78.33,4.11){\makebox(0,0)[cc]{$l(\pi_2,e_2)$}}
\put(104.66,3.78){\makebox(0,0)[cc]{$l(\pi_{n-1},e_{n-1})$}}
\put(139.00,3.11){\makebox(0,0)[cc]{$l(\pi_n,e_{n-1})$}}
\put(88.66,30.11){\makebox(0,0)[cc]{$p$}}
\put(88.66,13.11){\makebox(0,0)[cc]{$q$}}
\put(49.33,25.11){\makebox(0,0)[cc]{$e_1$}}
\put(33.33,21.11){\circle*{0.94}}
\put(46.66,21.11){\circle*{1.33}}
\put(60.33,21.11){\circle*{0.94}}
\put(129.66,20.78){\circle*{1.33}}
\put(73.66,21.11){\circle*{1.33}}
\put(117.00,20.78){\circle*{0.67}}
\put(95.00,20.56){\makebox(0,0)[cc]{...}}
\end{picture}
\begin{center}
\nopagebreak[4]
Fig. \theppp.

\end{center}
\addtocounter{ppp}{1}

The broken line formed of the intervals
$\ell(\pi_i,e_i)$, $\ell(\pi_i,e_{i-1})$ is called the {\em median} of this
band.

We say that two bands {\em cross} if their medians cross. We say that a band is
an {\em annulus} if its median is a closed curve.
In this case the first and the last cells of the band coincide (see \setcounter{pdeleven}{\value{ppp}}
Figure \thepdeleven)

\bigskip
\begin{center}
\unitlength=1.5mm
\linethickness{0.4pt}
\begin{picture}(101.44,22.89)
\put(30.78,13.78){\oval(25.33,8.44)[]}
\put(30.78,13.89){\oval(34.67,18.00)[]}
\put(39.67,9.56){\line(0,-1){4.67}}
\put(25.67,9.56){\line(0,-1){4.67}}
\put(18.56,14.89){\line(-1,0){5.11}}
\put(32.56,7.11){\makebox(0,0)[cc]{$\pi_1=\pi_n$}}
\put(21.00,7.33){\makebox(0,0)[cc]{$\pi_2$}}
\put(15.89,11.78){\makebox(0,0)[cc]{$\pi_3$}}
\put(43.44,12.89){\line(1,0){4.67}}
\put(43.00,8.44){\makebox(0,0)[cc]{$\pi_{n-1}$}}
\put(19.22,10.89){\line(-1,-1){4.00}}
\put(25.89,20.44){\circle*{0.00}}
\put(30.78,20.44){\circle*{0.00}}
\put(30.33,1.56){\makebox(0,0)[cc]{a}}
\put(35.44,20.44){\circle*{0.00}}
\put(62.78,4.89){\line(0,1){4.67}}
\put(62.78,9.56){\line(1,0){38.67}}
\put(101.44,9.56){\line(0,-1){4.67}}
\put(101.44,4.89){\line(-1,0){38.67}}
\put(82.11,16.22){\oval(38.67,13.33)[t]}
\put(62.78,15.78){\line(0,-1){7.33}}
\put(101.44,16.44){\line(0,-1){8.22}}
\put(82.00,12.22){\oval(27.78,14.67)[t]}
\put(67.89,12.89){\line(0,-1){8.00}}
\put(95.89,13.11){\line(0,-1){8.22}}
\put(67.89,13.78){\line(-1,1){4.67}}
\put(95.67,14.00){\line(1,1){5.11}}
\put(73.67,9.56){\line(0,-1){4.67}}
\put(89.67,9.56){\line(0,-1){4.67}}
\put(76.56,21.33){\circle*{0.00}}
\put(81.44,21.33){\circle*{0.00}}
\put(86.11,21.33){\circle*{0.00}}
\put(65.22,12.22){\makebox(0,0)[cc]{$\pi_1$}}
\put(65.22,6.89){\makebox(0,0)[cc]{$\pi$}}
\put(71.00,6.89){\makebox(0,0)[cc]{$\gamma_1$}}
\put(92.78,7.11){\makebox(0,0)[cc]{$\gamma_m$}}
\put(98.56,7.11){\makebox(0,0)[cc]{$\pi'$}}
\put(98.33,12.67){\makebox(0,0)[cc]{$\pi_n$}}
\put(94.11,20.22){\makebox(0,0)[cc]{$S$}}
\put(85.89,6.89){\makebox(0,0)[cc]{$T$}}
\put(76.56,6.89){\circle*{0.00}}
\put(79.44,6.89){\circle*{0.00}}
\put(82.33,6.89){\circle*{0.00}}
\put(83.44,1.56){\makebox(0,0)[cc]{b}}
\end{picture}
\end{center}
\begin{center}
\nopagebreak[4]
Fig. \theppp.

\end{center}
\addtocounter{ppp}{1}

Let $S$ and $T$ be two disjoint sets of letters, let ($\pi$, $\pi_1$, \ldots,
$\pi_n$, $\pi'$) be an $S$-band and let ($\pi$, $\gamma_1$, \ldots, $\gamma_m$,
$\pi'$) is a $T$-band. Suppose that:
\begin{itemize}
\item the medians of these bands form a simple closed
curve,
\item on the boundary  of $\pi$ and on the boundary of $\pi'$ the pairs of $S$-edges separate the pairs of $T$-edges,
\item the start and end edges of these bands are not contained in the
region bounded by the medians of the bands.
\end{itemize}
Then we say that these bands form an {\em $(S,T)$-annulus} and the curve formed by the medians of these bands is the {\em median}
of this annulus.

For example, the diagram on Figure \thepdsix\ contains $k$-bands, $q_i$-bands, $A$-bands crossing the circles transversally, and $r$-annuli filling the space between consecutive circles. In the diagram on Figure \thepdseven\, we also have a $\rho$-annulus going around the disc, and many $b$-bands consisting of the $b$-commutativity cells.

The main idea is  the following. In most relations of the presentation of $H$ one can choose
two pairs of letters which belong to disjoint sets of letters. For example, the relation $aba\iv b\iv=1$ has a pair of $a$-letters and a pair of $B$-letters. The cells corresponding to these
relations must form $a$-bands and $b$-bands in a van Kampen diagram. Each cell is an intersection of an $a$-band and a $b$-band. Thus if we prove that the number of $a$-bands is ``small" and the number of $b$-bands is ``small", and that an $a$-band and a $b$-band can have at most one common cell, we show that the number of $(a,b)$-commutativity cells is ``small".

In order to bound the number of bands we use the following idea. Suppose that we have ruled out annuli. Then every band starts (ends) either
on the boundary of the diagram (the number of such bands is linear in terms on the perimeter), or on the boundary of a cells (for example, an $a$-band can end on a disc). This gives us the direction in which to proceed.

First we assume that a diagram contains no discs and prove the absence of certain types of annuli: $r$-annuli, $\rho$-annuli, $a$-annuli, $(r,a)$-annuli, etc. (22 different kinds altogether). One way to prove it is to use a simultaneous induction: assume that one of these
annuli exists, take the innermost annulus of one of these kinds. Then the subdiagram bounded
by this annulus does not contain annuli of any of the 22 kinds. This makes the subdiagram look
nice and eventually leads to existence of a pair of cells that cancel (thus the diagram is not reduced which contradicts its minimality).

Then we assume that the diagram contains discs and we bound the number of discs (see below) and their perimeters. Then we bound the number of $r$-bands by proving that there are no $r$-annuli, so each of the $r$-bands must start and end on the boundary of the diagram. Similarly we bound the number of $\rho$-bands. Then we bound the
number of $q$-bands (they can start on the discs, and the perimeters of the discs are already bounded). Since every $q$-cell is an intersection of an $r$-band or a $\rho$-band and a $q$-band, we bound the number of $q$-cells. This leads to a bound of the number of $A$-bands (they can end on $q$-cells and on discs), and so on.

Of course we always need the absence of multiple intersections of bands. Although the next
\setcounter{pdtwelve}{\value{ppp}}
Figure \thepdtwelve\
shows that
a multiple intersection of an $S$-band
and a $T$-band  does not necessarily produce an $(S,T)$-annulus, it turns out to be enough to
rule out $(S,T)$-annuli.

\vskip -.5 in

\begin{center}
\unitlength=1.50mm
\linethickness{0.4pt}
\begin{picture}(97.67,32.89)
\put(11.67,10.67){\framebox(24.00,3.56)[cc]{}}
\bezier{88}(35.67,14.22)(38.78,24.67)(42.33,14.22)
\bezier{200}(42.33,14.22)(45.22,-2.67)(15.22,10.67)
\bezier{168}(32.33,14.22)(39.22,32.89)(46.11,12.22)
\bezier{240}(46.11,12.22)(48.33,-6.89)(11.67,10.67)
\put(15.22,10.67){\line(0,1){3.56}}
\put(32.33,14.22){\line(0,-1){3.56}}
\put(13.44,12.22){\makebox(0,0)[cc]{$\pi$}}
\put(34.11,12.45){\makebox(0,0)[cc]{$\pi'$}}
\put(24.11,12.45){\makebox(0,0)[cc]{$\ldots$}}
\put(36.78,17.33){\line(-1,1){2.40}}
\put(43.44,16.00){\circle*{0.00}}
\put(44.33,12.67){\circle*{0.00}}
\put(43.89,9.11){\circle*{0.00}}
\put(18.56,14.22){\line(0,-1){3.56}}
\put(68.78,10.67){\framebox(20.44,3.56)[cc]{}}
\put(72.33,10.67){\line(0,1){3.56}}
\put(72.33,14.22){\line(0,0){0.00}}
\put(75.67,14.22){\line(0,-1){3.56}}
\put(85.00,14.22){\line(0,-1){3.56}}
\put(70.56,12.44){\makebox(0,0)[cc]{$\pi$}}
\put(87.22,12.22){\makebox(0,0)[cc]{$\pi'$}}
\bezier{80}(89.22,14.22)(92.33,22.00)(93.00,10.67)
\bezier{80}(93.00,10.67)(94.11,6.44)(78.78,6.44)
\bezier{100}(84.33,6.67)(65.00,4.89)(65.22,10.67)
\bezier{72}(65.22,10.89)(65.22,20.89)(68.78,14.22)
\bezier{168}(85.00,14.22)(95.00,32.00)(96.11,10.67)
\bezier{108}(96.11,10.67)(97.67,2.89)(78.78,2.89)
\bezier{128}(83.89,3.11)(62.11,0.89)(62.33,10.67)
\bezier{156}(62.33,10.67)(63.89,30.89)(72.33,14.22)
\put(90.11,16.22){\line(-1,1){2.67}}
\put(66.33,16.89){\line(1,1){2.89}}
\put(80.56,12.22){\makebox(0,0)[cc]{$\ldots$}}
\put(81.00,4.67){\makebox(0,0)[cc]{$\dots$}}
\end{picture}
\end{center}
\begin{center}
\nopagebreak[4]
\vskip -.3 in
Fig. \theppp.

\end{center}
\addtocounter{ppp}{1}

In order to bound the number of discs (and their perimeters) in a \vk diagram, we use the following idea.

The generic diagram over the presentation of $H$ looks like this:

\bigskip
\unitlength=1mm
\linethickness{0.4pt}
\begin{picture}(112.33,67.00)
\put(51.33,20.00){\circle{13.50}}
\put(73.33,41.00){\circle{14.00}}
\put(87.67,17.67){\circle{14.00}}
\put(94.33,46.33){\circle{14.00}}
\put(47.67,46.00){\circle{14.00}}
\put(54.67,45.33){\line(6,-1){12.00}}
\put(66.67,43.33){\line(0,0){0.00}}
\put(53.67,43.00){\line(6,-1){12.67}}
\put(52.33,51.67){\line(1,4){4.00}}
\put(48.33,52.67){\line(1,4){3.67}}
\put(43.00,51.33){\line(-1,3){5.00}}
\put(41.33,49.00){\line(-1,3){5.33}}
\put(40.67,45.67){\line(-1,0){6.67}}
\put(41.33,43.00){\line(-1,0){7.33}}
\put(48.67,39.00){\line(1,-6){2.00}}
\put(46.00,39.00){\line(1,-6){2.33}}
\put(82.67,22.33){\line(0,0){0.00}}
\put(88.33,24.67){\line(1,3){5.00}}
\put(91.00,24.00){\line(1,3){5.33}}
\put(66.33,10.33){\circle{9.43}}
\put(79.67,43.67){\line(5,2){8.00}}
\put(80.33,41.33){\line(2,1){7.33}}
\put(71.67,48.00){\line(0,1){19.00}}
\put(74.33,48.00){\line(0,1){19.00}}
\put(57.00,16.33){\line(2,-1){5.67}}
\put(56.00,15.00){\line(2,-1){6.00}}
\put(71.00,11.67){\line(5,3){10.00}}
\put(71.00,9.67){\line(5,3){10.33}}
\put(87.00,10.67){\line(0,-1){6.67}}
\put(88.67,10.67){\line(0,-1){6.67}}
\put(67.67,25.33){\circle{8.67}}
\put(64.67,28.33){\line(-1,1){12.33}}
\put(51.00,39.67){\line(1,-1){12.67}}
\put(70.00,29.33){\line(1,3){1.67}}
\put(71.67,27.67){\line(2,5){2.67}}
\put(71.67,23.67){\line(5,-1){10.33}}
\put(71.00,22.33){\line(5,-1){10.33}}
\put(68.33,21.00){\line(0,-1){6.33}}
\put(67.00,21.00){\line(0,-1){6.33}}
\put(57.00,23.67){\line(6,1){6.33}}
\put(58.00,21.67){\line(5,1){6.33}}
\put(73.17,35.50){\oval(78.33,63.00)[]}
\put(92.67,53.33){\line(0,1){13.67}}
\put(95.00,53.33){\line(0,1){13.67}}
\put(72.00,25.67){\line(6,5){18.00}}
\put(72.00,24.33){\line(6,5){19.00}}
\put(68.67,46.00){\line(-1,5){4.33}}
\put(67.67,45.00){\line(-1,4){5.67}}
\put(77.00,47.00){\line(1,4){5.00}}
\put(78.33,46.00){\line(1,4){5.33}}
\put(101.33,48.00){\line(1,0){11.00}}
\put(101.33,45.33){\line(1,0){11.00}}
\put(98.33,52.33){\line(5,6){11.00}}
\put(100.00,50.67){\line(4,5){11.00}}
\put(94.67,19.33){\line(1,0){17.67}}
\put(94.67,17.00){\line(1,0){17.67}}
\put(92.67,13.00){\line(1,-1){9.33}}
\put(91.67,11.67){\line(1,-1){7.67}}
\put(62.67,7.33){\line(-3,-2){5.00}}
\put(63.67,6.67){\line(-5,-3){4.33}}
\put(65.33,5.67){\line(0,-1){1.67}}
\put(66.67,4.00){\line(0,1){1.67}}
\put(69.33,7.00){\line(5,-6){2.33}}
\put(73.00,4.00){\line(-3,4){2.67}}
\put(44.67,21.00){\line(-1,0){10.67}}
\put(34.00,19.67){\line(1,0){10.67}}
\put(45.67,16.00){\line(-4,-3){11.00}}
\put(36.00,6.33){\line(5,4){11.00}}
\put(50.00,13.33){\line(-2,-5){3.67}}
\put(47.67,4.00){\line(2,5){3.67}}
\end{picture}
\setcounter{pdthirteen}{\value{ppp}}
\begin{center}
\nopagebreak[4]
Fig. \theppp.

\end{center}
\addtocounter{ppp}{1}

Discs in the diagram are connected by $k$-bands.

So with every van Kampen diagram we can associate a graph of discs. The vertices of this graph are the discs plus one external vertex. Vertices are connected by the $k$-bands. If a $k$-band starts on a disc and ends on the boundary of the diagram, we assume that this band terminates in the external vertex. The degree of each internal vertex of this graph is $N>>1$. We prove that
this graph cannot have bigons: two discs connected by a pair of $k$-bands. This
implies that the graph of discs is {\em hyperbolic}, and a standard small cancellation theory applies \cite{LS}. In particular there exists a disc with $N-3$ external edges.  This also implies
that the number of discs and $k$-bands in the diagram is linear in terms of
the perimeter.

In order to rule out $r$-annuli, $\rho$-annuli and other types of annuli,
we use several type reducing surgeries on a diagram.
One of them is illustrated by the following picture.

 {\bf Moving $r$-bands.} Suppose that in a minimal diagram $\Delta$
an $r$-band $\rr$ touches a disc $\Pi$ as in
\setcounter{pdfourteen}{\value{ppp}}
Figure \thepdfourteen\ (that is one of the sides of $\rr$ has two common $k$-edges with the contour of the disc).
Then it can be proved that
the bottom path of $\rr$ has a common subpath with the contour of $\Pi$
starting and ending with $k$-edges. Let $p$ be the maximal common
subpath with this property, so that $\bott(R)=qpq'$, $\partial(\Pi)= p
p_1$.  Without loss of generality we can assume that the label $\Lab(p)$ of
the path $p$ has the form $k_iwk_{i+1}w'k_{i+2}...k_j$
where $w=uq_1uq_2uq_3$
Then for some word $V$ we have that
$\Lab(p)V$ is a cyclic shift of $K(w)=\Lab(\partial(\Pi))$. One can construct
an $r$-band $\rr'$ with the bottom path labeled by the word $V$
and the $r$-edges having the same labels as in $\rr$.  Let $\rr''$ be
the subband of $\rr$ with bottom path $p$, so $\rr=\rr_1\rr''\rr_2$. Let $e$
be the start edge of $\rr''$ and let $e'$ be the end edge of $\rr''$.  Cut
the diagram $\Delta$ along the path $e\iv p_1 e'$. We can fill the resulting
hole by gluing in the $r$-band $\rr'$ and the mirror image
$\overline\rr'$ of $\rr'$. The new diagram $\Delta'$ that we obtain this way
will have two $r$-bands instead of the old $r$-band $\rr$. The
first is $\rr_1(\overline\rr')\iv\rr_2$ (the inverse band
$(\overline\rr')\iv$ differs from $\overline\rr'$ by the order of  cells)
and the second one is $\rr''\rr'$. The second $r$-band is an annulus
which touches $\Pi$ along its inner boundary. If we replace the disc $\Pi$
by the corresponding \vk diagram over the presentation of $G_N(\sss)$, we
see that the subdiagram $\Pi'$ bounded by the outer boundary of the annulus
$\rr''\rr'$ is a diagram over the presentation of $G_N(\sss)$ with exactly
one hub and no $r$-edges on the boundary. Then one can prove that
$\Pi'$ is a disc (corresponding to some computation).
We replace it by one cell of the infinite presentation of $H$.
Then we reduce the resulting diagram.

\begin{center}
\unitlength=1.50mm
\linethickness{0.6pt}
\begin{picture}(84.50,56.00)
\put(7.33,33.08){\oval(9.33,15.83)[l]}
\put(20.16,33.08){\oval(15.33,15.83)[r]}
\put(7.33,41.00){\line(1,0){13.33}}
\bezier{64}(9.33,51.67)(7.66,43.50)(14.83,41.00)
\bezier{52}(10.66,51.67)(9.83,44.67)(15.83,42.50)
\bezier{60}(16.00,42.50)(25.16,43.33)(28.66,39.17)
\bezier{72}(28.83,39.17)(32.66,30.67)(26.16,24.67)
\bezier{56}(26.16,24.50)(24.33,22.50)(13.00,23.00)
\bezier{60}(13.00,22.83)(5.50,21.17)(5.83,13.83)
\bezier{72}(3.66,13.83)(3.33,22.67)(12.00,25.17)
\put(3.66,13.83){\line(1,0){2.17}}
\put(9.33,51.67){\line(1,0){1.33}}
\put(35.33,32.17){\vector(1,0){6.67}}
\put(59.00,32.58){\oval(9.33,15.83)[l]}
\put(71.83,32.58){\oval(15.33,15.83)[r]}
\put(59.00,40.50){\line(1,0){13.33}}
\bezier{60}(67.66,42.00)(76.83,42.83)(80.33,38.66)
\bezier{72}(80.50,38.66)(84.33,30.16)(77.83,24.16)
\bezier{56}(77.83,24.00)(76.00,22.00)(64.66,22.50)
\bezier{84}(67.16,44.16)(68.66,56.00)(77.33,54.50)
\bezier{84}(67.66,42.16)(69.83,54.83)(77.50,52.83)
\bezier{80}(63.66,20.33)(65.33,9.50)(74.33,9.83)
\bezier{76}(64.66,22.33)(65.66,12.50)(74.33,11.66)
\bezier{92}(67.66,42.33)(52.50,44.33)(52.16,36.16)
\bezier{80}(52.16,36.16)(50.66,23.00)(57.16,22.16)
\bezier{32}(57.00,22.16)(63.16,22.16)(64.66,22.33)
\bezier{132}(67.16,44.16)(48.66,47.16)(49.83,32.50)
\bezier{116}(50.00,32.50)(49.33,17.66)(63.66,20.33)
\put(64.66,32.33){\makebox(0,0)[cc]{Disc}}
\put(15.50,32.67){\makebox(0,0)[cc]{Disc}}
\put(25.16,33.00){\makebox(0,0)[cc]{$p$}}
\put(4.16,33.67){\makebox(0,0)[cc]{$p_1$}}
\put(56.16,32.83){\makebox(0,0)[cc]{$p_1$}}
\put(76.83,32.50){\makebox(0,0)[cc]{$p$}}
\put(6.83,46.00){\makebox(0,0)[cc]{$q$}}
\put(2.83,20.00){\makebox(0,0)[cc]{$q'$}}
\put(13.16,48.00){\makebox(0,0)[cc]{$\rr_1$}}
\put(32.83,33.33){\makebox(0,0)[cc]{$\rr''$}}
\put(8.83,16.33){\makebox(0,0)[cc]{$\rr_2$}}
\put(73.50,49.16){\makebox(0,0)[cc]{$\rr_1$}}
\put(84.50,32.16){\makebox(0,0)[cc]{$\rr''$}}
\put(70.83,16.16){\makebox(0,0)[cc]{$\rr_2$}}
\put(46.33,31.83){\makebox(0,0)[cc]{$\overline\rr'$}}
\put(14.33,41.00){\line(0,1){2.33}}
\put(7.33,25.17){\line(1,0){15.00}}
\put(12.16,25.17){\line(0,-1){2.50}}
\put(59.16,24.66){\line(1,0){13.33}}
\put(67.66,40.50){\line(0,1){1.67}}
\put(64.66,22.33){\line(0,1){2.33}}
\put(64.66,22.33){\rule{0.17\unitlength}{2.33\unitlength}}
\put(74.33,9.66){\rule{0.20\unitlength}{2.00\unitlength}}
\put(64.66,22.33){\line(-1,-2){1.20}}
\put(67.67,42.16){\line(-1,2){1.20}}
\put(67.83,42.16){\line(-1,2){1.20}}
\put(67.67,42.16){\line(-1,2){1.20}}
\put(77.50,52.66){\line(0,1){2.20}}
\end{picture}
\end{center}
\begin{center}
\nopagebreak[4]
\vskip -0.6 in
Fig. \theppp.

\end{center}
\addtocounter{ppp}{1}

This construction amounts
to changing the disc $\Pi$, moving the band $\rr$ through the disc and then reducing the resulting diagram.

Now suppose that there exists an $r$-annulus in our van Kampen diagram as in
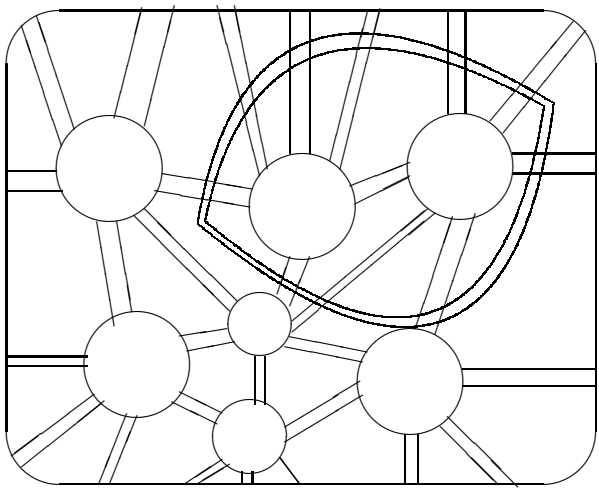
\setcounter{pdfifteen}{\value{ppp}}
Figure \thepdfifteen.

\unitlength=1.00mm
\linethickness{0.4pt}
\begin{picture}(112.34,79.33)
\put(51.33,20.00){\circle{13.50}}
\put(73.33,41.00){\circle{14.00}}
\put(87.67,17.67){\circle{14.00}}
\put(94.33,46.33){\circle{14.00}}
\put(47.67,46.00){\circle{14.00}}
\put(54.67,45.33){\line(6,-1){12.00}}
\put(66.67,43.33){\line(0,0){0.00}}
\put(53.67,43.00){\line(6,-1){12.67}}
\put(52.33,51.67){\line(1,4){4.00}}
\put(48.33,52.67){\line(1,4){3.67}}
\put(43.00,51.33){\line(-1,3){5.00}}
\put(41.33,49.00){\line(-1,3){5.33}}
\put(40.67,45.67){\line(-1,0){6.67}}
\put(41.33,43.00){\line(-1,0){7.33}}
\put(48.67,39.00){\line(1,-6){2.00}}
\put(46.00,39.00){\line(1,-6){2.33}}
\put(82.67,22.33){\line(0,0){0.00}}
\put(88.33,24.67){\line(1,3){5.00}}
\put(91.00,24.00){\line(1,3){5.33}}
\put(66.33,10.33){\circle{9.43}}
\put(79.67,43.67){\line(5,2){8.00}}
\put(80.33,41.33){\line(2,1){7.33}}
\put(71.67,48.00){\line(0,1){19.00}}
\put(74.33,48.00){\line(0,1){19.00}}
\put(57.00,16.33){\line(2,-1){5.67}}
\put(56.00,15.00){\line(2,-1){6.00}}
\put(71.00,11.67){\line(5,3){10.00}}
\put(71.00,9.67){\line(5,3){10.33}}
\put(87.00,10.67){\line(0,-1){6.67}}
\put(88.67,10.67){\line(0,-1){6.67}}
\put(67.67,25.33){\circle{8.67}}
\put(64.67,28.33){\line(-1,1){12.33}}
\put(51.00,39.67){\line(1,-1){12.67}}
\put(70.00,29.33){\line(1,3){1.67}}
\put(71.67,27.67){\line(2,5){2.67}}
\put(71.67,23.67){\line(5,-1){10.33}}
\put(71.00,22.33){\line(5,-1){10.33}}
\put(68.33,21.00){\line(0,-1){6.33}}
\put(67.00,21.00){\line(0,-1){6.33}}
\put(57.00,23.67){\line(6,1){6.33}}
\put(58.00,21.67){\line(5,1){6.33}}
\put(73.17,35.50){\oval(78.33,63.00)[]}
\put(92.67,53.33){\line(0,1){13.67}}
\put(95.00,53.33){\line(0,1){13.67}}
\put(72.00,25.67){\line(6,5){18.00}}
\put(72.00,24.33){\line(6,5){19.00}}
\put(68.67,46.00){\line(-1,5){4.33}}
\put(67.67,45.00){\line(-1,4){5.67}}
\put(77.00,47.00){\line(1,4){5.00}}
\put(78.33,46.00){\line(1,4){5.33}}
\put(101.33,48.00){\line(1,0){11.00}}
\put(101.33,45.33){\line(1,0){11.00}}
\put(98.33,52.33){\line(5,6){11.00}}
\put(100.00,50.67){\line(4,5){11.00}}
\put(94.67,19.33){\line(1,0){17.67}}
\put(94.67,17.00){\line(1,0){17.67}}
\put(92.67,13.00){\line(1,-1){9.33}}
\put(91.67,11.67){\line(1,-1){7.67}}
\put(62.67,7.33){\line(-3,-2){5.00}}
\put(63.67,6.67){\line(-5,-3){4.33}}
\put(65.33,5.67){\line(0,-1){1.67}}
\put(66.67,4.00){\line(0,1){1.67}}
\put(69.33,7.00){\line(5,-6){2.33}}
\put(73.00,4.00){\line(-3,4){2.67}}
\put(44.67,21.00){\line(-1,0){10.67}}
\put(34.00,19.67){\line(1,0){10.67}}
\put(45.67,16.00){\line(-4,-3){11.00}}
\put(36.00,6.33){\line(5,4){11.00}}
\put(50.00,13.33){\line(-2,-5){3.67}}
\put(47.67,4.00){\line(2,5){3.67}}
\bezier{320}(60.33,39.00)(67.67,75.33)(105.33,54.33)
\bezier{384}(60.33,39.00)(98.00,7.33)(105.33,54.00)
\bezier{356}(59.33,38.67)(66.00,79.33)(106.67,54.67)
\bezier{416}(59.33,38.67)(100.33,4.67)(106.67,54.67)
\end{picture}
\begin{center}
\nopagebreak[4]
Fig. \theppp.

\end{center}
\addtocounter{ppp}{1}

Suppose that there are discs inside the region bounded by this annulus. Then
these discs form a hyperbolic graph, and so
the $r$-band will intersect more than 1/2 of the $k$-bands going out of
one of these discs. Then the $r$-band moving construction reduces the type of the
diagram. Thus the region bounded by the $r$-annulus cannot contain discs. But
we have ruled out the case when a diagram without discs contains an $r$-annulus, a contradiction.

In order to bound the perimeters of discs and $B$-cells we use the following
idea. The contour of a disc contains a constant number of non $A$-edges. Thus in order to
bound the perimeter of a disc, we need to bound the number of $A$-edges on the contour of it.
Every $A$-edge on the contour of a discs is the start edge of an $a$-band. An $a$-band consists
of $a$-commutativity cells corresponding to relations of the form $ab=ba$, $ar=ra$ or $a\rho=\rho a$ or to the relations of the form $a^\rho=ab$. Thus an $a$-band can end either on a disc or on the boundary of a $(a,q,r)$-cell. The latter belongs to an $r$-band and we already know that the diagram contains only ``small number" of $r$-bands. Thus if the a disc has a very big perimeter and many of the $a$-bands starting on the contour of this disc end on  boundaries of $(a,q,r)$-cells, then many of these $a$-bands must end on the contour of the same $r$-band. The following lemma shows that it is impossible.

\begin{lemma} Let $\rr_1$,...,$\rr_n$ be maximal $a$-bands starting
on a path $p$ where $p$ is an $A$-subpath of the boundary of
a disc $\Pi$. Suppose that the end edges of all $\rr_i$ are on
the contours of $r$-cells belonging to the same
$r$-band $\ttt$.
Then $n\le 2$.
\label{lm}
\end{lemma}

\noindent {\bf Sketch of the Proof.} Indeed, if $n>2$
then there are three $a$-bands, say,  $\rr_1$, $\rr_2$, $\rr_3$
starting on $p$ and ending on three different cells $\pi_1$, $\pi_2$ and
$\pi_3$ of $\ttt$. We can assume that $\pi_2$ is between $\pi_1$ and
$\pi_3$ (see
\setcounter{pdsixteen}{\value{ppp}}
Figure \thepdsixteen).  Consider the minimal subdiagram $\Delta_1$ of our diagram
containing  $a$-bands $\rr_1, \rr_2, \rr_3$, the minimal subpath
of the path $p$ containing the starting edges of $\rr_1$, $\rr_2$, $\rr_3$,
 and the part of the band $\ttt$ between $\pi_1$ and $\pi_3$
Then $\Delta_1$ has
 no $k$-edges on its contour. Therefore $\Delta_1$ does not contain
 discs. Therefore the maximal
 $q$-band $\qq$ in $\Delta_1$ containing $\pi_2$ divides
$\Delta_2$ into two parts (that is if we delete the $q$-edges
 from $\qq$, the diagram $\Delta_2$ will fall into two pieces).  The subpath
 of the path $p$ containing the start edges of $\rr_1,\rr_2,\rr_3$ is
 contained in one of these parts since it does not contain $\overline
 Q$-edges.  The cells $\pi_1$ and $\pi_3$ belong to different parts because
 $\qq$ cannot intersect $\ttt$ twice. Since the
 $\pi_1$ and $\pi_3$ are connected with the cells on $p$ by $a$-bands, one of these bands must intersect $\qq$ which is impossible (a
 $q$-band cannot cross an $a$-band).  $\Box$

\medskip
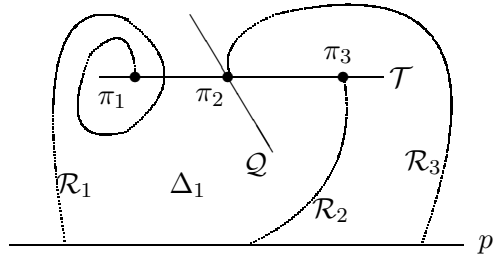
\begin{center}
\unitlength=1.50mm
\linethickness{0.4pt}
\begin{picture}(45.67,25.78)
\put(0.11,0.22){\line(1,0){40.22}}
\put(8.11,15.11){\line(1,0){25.11}}
\put(11.22,15.11){\circle*{0.99}}
\put(29.67,15.11){\circle*{0.89}}
\bezier{136}(5.00,0.22)(1.22,23.56)(10.56,20.22)
\bezier{68}(10.33,20.22)(17.22,15.33)(10.33,10.22)
\bezier{64}(10.33,10.22)(3.89,8.67)(7.22,17.11)
\bezier{40}(7.22,17.11)(10.78,20.67)(11.22,15.33)
\bezier{200}(36.56,0.22)(45.67,25.78)(23.67,20.67)
\bezier{36}(23.67,20.67)(18.33,18.44)(19.67,15.11)
\put(19.44,15.11){\circle*{0.89}}
\bezier{84}(21.22,0.22)(31.44,5.11)(29.67,15.11)
\put(16.11,20.67){\line(3,-5){7.33}}
\put(9.22,13.11){\makebox(0,0)[cc]{$\pi_1$}}
\put(17.89,13.33){\makebox(0,0)[cc]{$\pi_2$}}
\put(29.22,17.11){\makebox(0,0)[cc]{$\pi_3$}}
\put(21.89,7.33){\makebox(0,0)[cc]{$\qq$}}
\put(34.78,15.11){\makebox(0,0)[cc]{$\ttt$}}
\put(5.89,5.56){\makebox(0,0)[cc]{$\rr_1$}}
\put(15.89,5.56){\makebox(0,0)[cc]{$\Delta_1$}}
\put(28.56,3.33){\makebox(0,0)[cc]{$\rr_2$}}
\put(36.78,7.56){\makebox(0,0)[cc]{$\rr_3$}}
\put(42.33,0.22){\makebox(0,0)[cc]{$p$}}
\end{picture}
\end{center}

\begin{center}
\nopagebreak[4]
Fig. \theppp.

\end{center}
\addtocounter{ppp}{1}

Finally we need to estimate the number and perimeters of $B$-cells (i.e. cells
corresponding to relations of the group $\langle B\rangle$). Here we use the following
trick. Suppose that two $B$-cells are connected by a $b$-band consisting of $(a,b)$-commutativity cell. Then we can cut the two $B$-cells together with the $b$-band from the
diagram, and  replace it by one $B$-cell and a number of $(a,b)$-commutativity relations.
This reduces the type of the diagram because the commutativity relations have smaller rank
than $B$-relations.
\setcounter{pdseventeen}{\value{ppp}}
Figure \thepdseventeen\ shows how this surgery proceeds.

\unitlength=1.50mm
\linethickness{0.4pt}
\begin{picture}(76.83,21.17)
\bezier{96}(24.17,11.17)(16.00,2.33)(7.33,11.17)
\bezier{92}(7.33,11.17)(14.83,19.83)(24.17,13.67)
\put(24.17,11.00){\line(0,1){2.33}}
\put(24.17,13.67){\line(1,0){21.50}}
\put(45.67,13.67){\line(0,-1){2.33}}
\put(45.67,11.33){\line(-1,0){21.50}}
\put(24.17,11.33){\line(0,1){2.33}}
\bezier{100}(45.67,13.67)(56.17,21.17)(64.83,12.67)
\bezier{100}(45.67,11.33)(56.83,3.83)(64.67,12.67)
\put(26.67,11.33){\line(0,1){2.33}}
\put(29.67,11.33){\line(0,1){2.33}}
\put(32.67,11.33){\line(0,1){2.33}}
\put(35.83,11.33){\line(0,1){2.33}}
\put(39.33,11.33){\line(0,1){2.33}}
\put(42.17,11.33){\line(0,1){2.33}}
\put(69.00,12.67){\vector(1,0){7.83}}
\put(22.83,12.17){\makebox(0,0)[cc]{$b$}}
\put(47.00,12.50){\makebox(0,0)[cc]{$b$}}
\put(55.50,15.50){\makebox(0,0)[cc]{$v_b$}}
\put(14.17,14.17){\makebox(0,0)[cc]{$u_b$}}
\put(34.33,15.17){\makebox(0,0)[cc]{$w_a$}}
\put(34.17,9.50){\makebox(0,0)[cc]{$w_a$}}
\end{picture}

\linethickness{0.4pt}
\begin{picture}(66.17,24.83)
\bezier{96}(36.00,11.33)(27.83,2.50)(19.17,11.33)
\bezier{92}(19.17,11.33)(26.67,20.00)(36.00,13.83)
\bezier{100}(36.00,13.67)(46.50,21.17)(55.17,12.67)
\bezier{100}(36.00,11.33)(47.17,3.83)(55.00,12.67)
\put(36.00,11.17){\line(0,1){2.50}}
\bezier{152}(36.00,1.67)(8.67,1.50)(8.33,12.33)
\bezier{140}(8.17,12.33)(14.67,24.83)(36.00,24.33)
\bezier{152}(35.83,24.33)(63.67,23.00)(66.17,12.67)
\bezier{160}(36.00,1.67)(64.50,1.50)(66.17,12.50)
\put(35.00,4.00){\makebox(0,0)[cc]{$(A,B)$-commutation cells}}
\put(35.17,12.50){\makebox(0,0)[cc]{$b$}}
\put(24.00,13.67){\makebox(0,0)[cc]{$u_b$}}
\put(45.33,15.33){\makebox(0,0)[cc]{$v_b$}}
\put(30.17,22.17){\makebox(0,0)[cc]{$u_bw_av_bw_a\iv$}}
\end{picture}

\begin{center}
\nopagebreak[4]
Fig. \theppp.

\end{center}
\addtocounter{ppp}{1}

This implies that every $b$-band starting on the contour of a $B$-cell must
end either on the boundary of the diagram or on the contour of a
$(a,\rho,b)$-cell. The number of maximal $a$-bands in the diagram is
bounded (because the total perimeter of the discs is bounded, and the number
of $q$-cells is bounded too), and a lemma similar to Lemma \ref{lm} shows
that the number of $b$-bands starting on the contour of the same $B$-cell
and ending on the contour of the same $a$-band is at most 2. This leads to
the bound of the number of $B$-cells and the total perimeter of $B$-cells.

Finally we can estimate the areas of words in $H$ relative to the finite
presentation of $H$. Take any word $w$ which is equal to 1 in $H$. Then
there exists a diagram over the infinite presentation of $H$ (with discs and
$B$-cells) with boundary label $w$. The total perimeter of discs and
$B$-cells is bounded by a polynomial in $|w|$. Now replace every disc by the
\vk diagram over the finite presentation of $H$ (as in Fig. \thepdsix), and
replace each $B$-cell by the diagram on Fig. \thepdseven\ consisting of two
discs and a relatively small number of other cells. The resulting diagram
will be a \vk diagram over the finite presentation of $H$. It is easy to see
that if the perimeter of a disc is $p$ then the area is $O(T(p)^2)$ where
$T$ is the time function of the $S$-machine. This gives an estimate of the
area of $w$ which is polynomially equivalent to $T(|w|)$.

\subsection{Why Is There No Distortion?}

The proof that the embedding of $B(m,n)$ into $H_{m,n}$ and in general any
recursively presented group $G$ into $H_N(\sss)$ is undistorted also uses
bands and annuli, and the structure of diagrams over the infinite
presentation of $H$ described in the previous section.

Here we present the main points of the proof of bounded distortion in
Theorems \ref{ol1}, \ref{th1} and \ref{sss}.

For simplicity consider the case of the group $H=H_{m,n}$ from Theorem
\ref{sss}. The general case of $H_N(\sss)$ is similar.  By definition of
bounded distortion, we have to find a constant $c>0$ such that for any
element $g\in B(m,n)$ represented by a geodesic (in $B(m,n)$) word
$U=U(b_1,\dots,b_m)$ in the alphabet $B=\{b_1,\dots,b_m\}$ and for any word
$Z$ in the generators of $H$, that represents the same element, we have
$|Z|\ge c|U|$.

In order to achieve this goal, consider the minimal diagram $\Delta$ over
the infinite presentation of $H$ considered in the previous section, with
boundary label $UZ\iv$.  Then the boundary of $\Delta$ has the form $p\iv
p'$ where $\Lab(p)\equiv U$, $\Lab(p')\equiv Z$.  We need to show that
$|p|\le c|p'|$ for some constant $c$. It suffices to make a correspondence
between $b$-edges of $p$ and edges of $p'$, such that any edge of $p'$
corresponds to at most $c\iv$ edges of $p$.

First of all notice that we can assume that no $B$-cell in $\Delta$ has a
common edge with $p$.  Indeed, if such a $B$-cell exists, we can cut it off
reducing the type of the diagram and replacing the path $p$ with a not
shorter path $p_1$ (recall that $U$ was a geodesic word representing $g$).

Therefore for every edge $e$ on $p$,  there is a maximal $b$-band ${\cal T}$
in $\Delta$, starting at $e$. It can end neither on $p$ nor on the boundary
of a $G_b$-cell (both cases are ruled out in the same manner as it was done
in the previous section: we can do a type reducing surgery again).

If ${\cal T}$ ends on the path $p'$, we associate the terminal edge of
$\ttt$ with $e$.  Another possibility is that ${\cal T}$ terminates on the
boundary of some maximal $a$-band (at the cell labeled by a relator
$\rho^{-1} a_i\rho b_i^{-1}a_i^{-1}$). A lemma similar to Lemma \ref{lm}
shows that at most 2 maximal $b$-bands starting on $p$ can end on the
boundary of the same $a$-band. This means that we can consider the set
$\aaa$ of $a$-bands where these $b$-bands end.

The most pleasant (for us) among these $a$-bands are those which start or
end on $p'$ (they cannot end on $p$ because $p$ does not have $A$-edges).
Other $a$-bands can terminate either on contours of $r$-bands or on disks.

We need to consider two cases. In the first case the number of those
$r$-bands is large (proportional to the number of $a$-bands in $\aaa$).
Since there are no $r$-annuli in $\Delta$ (see the previous section), each
of these $r$-bands starts and ends on $p'$, and we obtain a desired
inequality $|p'|\ge c|p|$.

In the second case we have to assume that the number of $r$-bands where
$a$-bands terminate is ``small". Since by a variation of Lemma \ref{lm} the
number of $a$-bands terminating on the same $r$-band is bounded by a
constant, in this case most $a$-bands in $\aaa$ terminate on discs.

In this situation we use the so called {\em ovals} and their {\em shadows}
(see \cite{Ol1} and \cite{BORS}).

An {\it oval} is a simple closed path $h$ in the disk graph of the diagram
$\Delta$. It divides the plane into two regions. One of them, denoted by
$O(h)$, must possess the following property. For every disk $\Pi$ on $h$,
the number $n_1$ of maximal $k$-bands going from $\Pi$ into $O(h)$ and the
number $n_2$ of the maximal $k$-bands going from $\Pi$ into the exterior of
$O(h)$ satisfy the following inequalities:  $$n_1\ge n_2+8, \quad n_2\le
2n_1.$$

The hyperbolicity of the disk graph and high degrees of its interior
vertices make possible drawing an oval $h$ passing via any interior edge of
the disc graph of $\Delta$.

One of the  main properties of ovals is that if an $r$-band starts outside
the subdiagram $O(h)$ bounded by the oval and  then intersects the oval,
it cannot leave $O(h)$. Indeed otherwise the hyperbolicity of the disk
graph would
imply the existence of either a  $(k,r)$-annulus (which is impossible, see
the previous section) or an $r$-band intersecting too many
maximal $k$-bands starting on the same disk (again it is impossible because
of the Moving $r$-bands construction from the previous section).

Thus any maximal $r$-band crossing an oval $h$, must intersect its {\it
shadow}, i.e. the boundary subpath of the diagram lying in $O(h)$.  This
allows us to prove that the shadow of any oval $h$ is sufficiently long
comparing to the perimeter of any disk $\Pi$ crossed by $h$. We can also
choose $h$ in such a way that the shadow of $h$ is inside $p'$ (because $p$
does not contain $k$-edges). Therefore the number of $a$-bands ending on the
contour of a disc does not exceed the length of the shadow of any oval
passing through the disc. If the bands $\aaa$ end on different discs
$\Pi_1$,...,$\Pi_k$ then the hyperbolicity of the disc graph allows us to
find ovals passing through these discs which have disjoint shadows, all
inside $p'$. Thus the length of $p'$ cannot be smaller than the number of
$a$-bands in $\aaa$, which in turn, as we know, cannot be much smaller than
the length of $p$.

The proof of the result that the shadow of an oval $h$ is sufficiently long
comparing to the  perimeter of a disc $\Pi$ in $h$ consists of two cases. In
the first case the number of maximal $r$-band in $O(h)$ is sufficiently
large (greater than, say,  $\frac{1}{20}$ of the number of all $a$-edges
between successive $k$-edges of $\Pi$). This case is clear since all the
$r$-bands must terminate on the shadow.

The second case is complementary to the first one. Since the number
of the $r$-bands is small, the quantity of the $a$-bands going from $\Pi$ into $O(p)$
and terminating on $r$-bands, is small too (Lemma \ref{lm} works again).
Therefore a majority of them terminates either in the shadow of $\Pi$
(this is the best alternative for us) or on some disks
$\Pi',\Pi'',\dots \Pi^{(k)}$ inside of $O(h)$.

This situation can be analyzed by induction: as before we can draw ovals
$h_1,\dots,h_k$ passing through $\Pi',\dots,\Pi^{(k)}$ respectively, whose shadows are disjoint and are inside the shadow of $h$.

\subsection{Embeddings With Given Length Functions}

Here we present the main ideas of the proofs of results from \cite{Ol4} (see
Section \ref{lf}).

If $G$ is a subgroup of a group $H$ with a finite set of generators
$B=\{b_1,\dots,b_n\}$ then the function $\ell(g)=|g|_B$ on $G$ evidently
satisfies conditions (D1)-(D3) from Theorem \ref{ol4}.
For instance condition (D3) holds because the
number of all words of length at most $k$ in the alphabet $B$ grows
exponentially as $k\rightarrow\infty$.

To prove that every function $G\to {\bf N}$ satisfying conditions (D1)-(D3)
can be realized as the length function of $G$ inside a finitely generated
group $H$, we start with a presentation $G = F_G/N$, where $F_G$ is a free
group with the basis $\{x_g\}_{g\in G\backslash \{1\}}$ and $N$ is the
kernel of homomorphism $\epsilon: x_g\mapsto g$.

Then we construct an embedding $\beta: x_g\to X_g$ of $F_G$ into the
2-generated free group $F=F_2=F(b_1,b_2)$, such that the image $\beta(F_G)$
is freely generated by $X_g$, $g \in G\backslash \{1\}$ and the words $X_g$
are very ``independent" in the sense described below.

The group $H$ is equal to the quotient $F/L$, where $L$ is the normal
closure of $\beta(N)$ in $F$.

Notice that whatever homomorphism $\beta$ we choose, it induces a
homomorphism $\gamma: G\rightarrow H$. To make this homomorphism injective,
we need the following property:  $$L\cap\beta(F_G) = \beta(N).$$ In fact
$\beta(F_G)$ satisfies the following much stronger property:

\medskip
(*) {\it For any normal subgroup $U\triangleleft \beta(F_G)$ there is a
normal subgroup $V\triangleleft F$ such that $U=V\cap \beta(F_G)$.}

\medskip

The fact that free groups and more generally every non-elementary hyperbolic
group has plenty of infinitely generated free subgroups with property (*) is
interesting in its own right, it was the key ingredient in Olshanskii's
proof from \cite{Ol6} of the fact that every non-elementary hyperbolic group
is SQ-universal.

It turns out that we can make $\beta$ satisfy condition (*) by choosing
reduced words $X_g$ with the following small cancellation condition:

\medskip
(**) {\it If $Y$ is a subword of a word $X_g$ and $|Y|\ge\frac{1}{50}|X_g|$
then $Y$ occurs in $X_g$  as a subword only once, and $Y$  occurs neither
in $X_g^{-1}$ nor in $X_h^{\pm 1}$ for $h\ne g$.}
\medskip

It is relatively easy to construct an infinite set of words
$X_g$ in the alphabet $\{b_1,b_2\}$ which satisfies
the (**)-condition and has exponential growth, that is the number of different words $X_g$
of length $k$ grows exponentially as $k\to \infty$.

Since by condition (D3) the number of elements $g\in G$ with $\ell(g)\le k$ does not exceed
$c^k$ for some constant $c$, we can choose
the set $\{X_g\}$ in such a way that
$$\ell(g)\le |X_g|<d\ell(g)$$ for some positive constant $d$
and every $g\in G\backslash \{1\}$.

We need to show that the embedding $\gamma: G\to H$ has bounded distortion. For this we take any element $X_g$ of $\gamma(G)$ and consider the shortest word $W$ in the alphabet $\{b_1, b_2\}$
representing $X_g$ in $H$. The group $H$ is given by the presentation consisting of all relations of the form $X_{g_1}X_{g_2}...X_{g_n}=1$ where $g_1g_2...g_n=1$ in $G$.

Since $X_g=W$ modulo this presentation, we can consider the corresponding \vk diagram $\Delta$
with boundary label $X_gW\iv$. We can assume that the number of cells in $\Delta$ is minimal
among all such diagrams.

The condition (**) implies the following property of \vk diagrams over the
presentation of $H$.  Let $\Pi_1$ and $\Pi_2$ be  cells in a diagram $\Delta$
having a common edge.  Then either any common arc $p$ of the boundaries
$\partial\Pi_1$ and $\partial\Pi_2$ is short comparing to the perimeters
$P_1, P_2$ of the cells (say, $|p|\le\frac{1}{10}\min(P_1, P_2)$), or a
subdiagram consisting of $\Pi_1$ and $\Pi_2$, has also a boundary label of
the form $w=w(X_{g_1},\dots,X_{g_n})$, i.e. the subdiagram can be replaced
by one cell. The latter option cannot occur in $\Delta$ because of the
minimality of the diagram $\Delta$. Thus $\Delta$ satisfies a small
cancellation condition \cite{LS}.

This in turn allows us to prove that the word $W$ is freely equal to a
product $X_{g_1}^{\pm 1}\dots X_{g_s}^{\pm 1}$ for some $g_j\in G$ with
$g=g_1^{\pm 1}\dots g_s^{\pm 1}$ (see \cite{Ol4} for details). Further,
since the cancellations in such a product are small, $$|\gamma(g)|_H\ge
(1-\frac{2}{50})\sum_{j=1}^s|X_{g_1}|.$$ By conditions (D1), (D2), and by
the choice of $X_g$, we have:  $$|\gamma(g)|_H\ge 0.96\sum
\ell(g_j)=0.96\sum \ell(g_j^{\pm 1})\ge 0.96 \ell(g).$$ Hence $0.96
\ell(g)\le |\gamma(g)|_H\le d\ell(g)$, so $\ell$ is $O$-equivalent to the
length function of $G$ in $H$.  \label{efg}

\bigskip
\noindent Alexander Yu. Olshanskii\\
Department of Higher Algebra\\
MEHMAT. Moscow State University\\
olsh@nw.math.msu.su\\

\noindent Mark V. Sapir\\
Department of Mathematics\\
Vanderbilt University\\
http://www.math.vanderbilt.edu/$\sim$msapir\\

\end{document}